\documentclass[11pt]{article}

\usepackage{amsmath}
\usepackage{amssymb}
\usepackage{latexsym}
\usepackage{amsmath, amsfonts,amssymb, amsthm, euscript,makeidx,color,mathrsfs}
\usepackage{xcolor}
\usepackage{indentfirst}
\usepackage{algpseudocode}
\usepackage{graphicx}
\usepackage{epstopdf}
\usepackage{float}
\usepackage{subfig}

\newtheorem{assumption}{Assumption}
\def\qed{ \ \vrule width.2cm height.2cm depth0cm\smallskip}

\newcommand{\ba}{\begin{array}}
\newcommand{\ea}{\end{array}}
\newcommand{\be}{\begin{equation}}
\newcommand{\ee}{\end{equation}}
\newcommand{\bea}{\begin{eqnarray}}
\newcommand{\eea}{\end{eqnarray}}
\newcommand{\beaa}{\begin{eqnarray*}}
\newcommand{\eeaa}{\end{eqnarray*}}

\def\neg{\negthinspace}

\def\a{\alpha}

\def\e{\varepsilon}

\def\l{\lambda}
\def\m{\mu}

\def\si{\sigma}

\def\th{\theta}
\def\o{\omega}

\def\D{\Delta}

\def\O{\Omega}

\def\D{\Delta}

\def\O{\Omega}
\def\cF{{\cal F}}
\def\cG{{\cal G}}

\def\cI{{\cal I}}

\def\cL{{\cal L}}

\def\cN{{\cal N}}

\def\hE{\mathbb{E}}
\def\hF{\mathbb{F}}

\def\hL{\mathbb{L}}

\def\hN{\mathbb{N}}

\def\hP{\mathbb{P}}
\def\hQ{\mathbb{Q}}
\def\hR{\mathbb{R}}

\def\hX{\mathbb{X}}

\def\sA{\mathscr{A}}
\def\sB{\mathscr{B}}
\def\sC{\mathscr{C}}

\def\sG{\mathscr{G}}

\def\sK{\mathscr{K}}

\def\sM{\mathscr{M}}

\def\sP{\mathscr{P}}

\def\sX{\mathscr{X}}

\def\no{\noindent}

\def\ss{\smallskip}
\def\ms{\medskip}

\def\q{\quad}
\def\qq{\qquad}

\def\cd{\cdot}
\def\cds{\cdots}

\def\bx{{\bf x}}

\def\qed{ \hfill \vrule width.25cm height.25cm depth0cm\smallskip}

\newcommand{\basa}{\begin{assumption}}
\newcommand{\easa}{\end{assumption}}

\newcommand{\ol}{\overline}

\newcommand{\bas}{\begin{assum}}
\newcommand{\eas}{\end{assum}}

\def\liminf{\mathop{\underline{\rm lim}}}

 \def\cd{\cdot}
\def\cds{\cdots}
\def\ae{\hbox{\rm-a.e.{ }}}
\def\as{\hbox{\rm-a.s.{ }}}

\def\co{\mathop{{\rm co}}}

\def\ol{\overline}

\def\cl{\mbox{\rm cl}}
\def\co{\mbox{\rm co}}

\def\dis{\displaystyle}

\def\bx{{\bf x}}

\def\1{{\bf 1}}

\def\:{\!:\!}

\DeclareMathOperator{\dec}{dec}

at 9pt

\def\prehp(#1,#2){\ensuremath{  #1 \cdot #2 }}

\begin{document}

\newtheorem{thm}{Theorem}[section]
\newtheorem{lem}[thm]{Lemma}
\newtheorem{cor}[thm]{Corollary}
\newtheorem{prop}[thm]{Proposition}
\newtheorem{rem}[thm]{Remark}
\newtheorem{eg}[thm]{Example}
\newtheorem{defn}[thm]{Definition}
\newtheorem{assum}[thm]{Assumption}

\renewcommand {\theequation}{\arabic{section}.\arabic{equation}}
\def\thesection{\arabic{section}}
\title{\bf Set-Valued Stochastic Differential Equations with Unbounded Coefficients}

\author{
Atiqah Almuzaini\thanks{\noindent Department of
Mathematics, University of Southern California, Los Angeles, 90089; email: ahalmuza@usc.edu}
~ and ~ Jin Ma\thanks{ \noindent Department of
Mathematics, University of Southern California, Los Angeles, 90089;
email: jinma@usc.edu. This author is supported in part
by US NSF grant  \#2205972. } }

\date{\today}
\maketitle

\date{}
\maketitle

\begin{abstract}
 In this paper, we mainly focus on the set-valued (stochastic) analysis on the space of convex, closed, but possibly unbounded sets, and try to establish a useful theoretical framework for studying the set-valued stochastic differential equations with unbounded coefficients. The space that we will be focusing on are convex, closed sets that are ``generated" by a given cone, in the sense that the Hausdorff distance of all elements to the ``generating" cone is finite. Such space should in particular include the so-called ``upper sets", and has many useful cases in finance, such as the well-known set-valued risk measures, as well as the solvency cone in some super-hedging problems. We shall argue that, for such a special class of unbounded sets, under some conditions,  the cancellation law is still valid, eliminating a major obstacle for extending the set-valued analysis to non-compact sets. We shall establish some basic algebraic and topological properties of such spaces, and show that some standard techniques will again be valid in studying the set-valued SDEs with unbounded (drift) coefficients which, to the best of our knowledge, is new.

\end{abstract}

\vfill
{\bf Keywords:}
Set-valued stochastic differential equations, unbounded coefficients, set-valued stochastic analysis on unbounded sets, $L_C$-space,  solvency cones.

{\bf AMS subject classifications:} 60H10, 28B20, 47H04,  52A07, 54C60, 91G80.

\newpage

\section{Introduction}

Set-valued analysis, both deterministic and stochastic, has been widely applied in 
optimizations, optimal control, as well as economics and finance. However,   the existing literature in set-valued analysis predominantly focuses on compact sets. Although such a limitation is often technically necessary and does not pose fundamental challenge in many applications, it becomes more pronounced when the underlying objects under investigation are unbounded in nature. 

There are many practical applications that involve set-valued subjects that are necessarily unbounded. A well-known example is the  {\it set-valued risk measures},  the extension of the univariate risk measure when the risk appears in the form of  random vectors, 
which is often referred to as {\it systemic risk} in the context of default contagion (see, e.g., \cite{conical}) or multi-asset markets with transaction costs (see, e.g,   \cite{rmbook},\cite{FRsurvey}).
Mathematically, the set-valued risk measure  takes the form of functions with values being the ``upper sets", that is, the convex sets that are  ``additively invariant" with respect to a  fixed convex cone, 
and hence unbounded in nature. It is worth noting that a univariate  {\it dynamic risk measure} (the  family of risk measures indexed by time) 
 satisfying the so-called {\it time-consistency}   is often representable by the {\it backward stochastic differential equation} (BSDE)  (cf. \cite{bionnadal,coquet,gianin}), but the extension of such representation seems to be quite remote without appropriate technical tools in set-valued analysis that can deal with unbounded sets.  Another  example that involves upper sets is the so-called ``solvency cone" proposed in, e.g.,  Kabanov  \cite{Kcurrency,KandS} and Scharchmayer \cite{Sch}, which are frequently used in the study of super-hedging problem that involves the transection costs.  It is quite conceivable that any dynamics 
 whose coefficients involve the solvency cone would naturally require a framework for unbounded sets. 
%
%

This paper is an attempt to establish a workable theoretical basis of set-valued analysis that is suitable for the study of set-valued SDEs that  have unbounded coefficients. There are, however, several main technical obstacles that need to be recognized, before we try to build the reasonable framework. 
The first obstacle regarding the  space of unbounded sets is the lack of the so-called cancellation law in its algebraic structure under the  Minkowski addition. More precisely,  for nonempty closed convex sets $A$, $B$, and $D\subseteq \mathbb{R}^d$ with $A+B=D+B,$ it is not true in general  that $A=D$, unless  $B$ is compact. The second main obstacle in
dealing with unbounded sets is in its topological structure with respect to the Hausdorff metric. More precisely, there is no single set (bounded or unbounded) that has a finite distance from all unbounded sets. In particular, the Hausdorff distance between  a given set and  any other (unbounded) set is equal to  infinity in general. Therefore, our first task is to try to identify an appropriate space of unbounded sets that possesses a  ``reference point", that is, an unbounded set   that has finite distance to every other elements in the  space. 

%
%

It turns out that the space of unbounded sets that have the ``upper-set" nature with respect to a certain fixed cone is one that fits our purpose well. Not only does such a space (which we call the $L_C$-space in this paper) contain many well-known application, we can actually argue that such a space of unbounded sets that can be equipped with appropriate topological and algebraic structures on which our set-valued analysis will be based, and that many desired properties of set-valued analysis can be established with some efforts. Moreover,  we shall extend some necessary concepts, such as set-valued Lebesgue integral and stochastic integrals and their path regularities under the new framework of $L_C$-valued unbounded sets. 
%

Another main objective of this paper is to establish the existence and uniqueness of the solution to a class of set-valued stochastic differential equations with bounded drift coefficients. We recall that the first papers dealing with deterministic set-valued  differential equations with compact set-valued coefficients can be traced back to the early works of \cite{AA95,DIF69,DBP70}. The stochastic set-valued  differential equations with compact  drift and absolutely summable diffusion terms and their path-regularity can be found in, e.g., \cite{KisMic16}. However, to the best of our knowledge, the set-valued SDEs or even stochastic differential inclusions (SDIs) with unbounded coefficients, even only on the drift part, is novel in the literature. We should remark that in this paper, we shall content ourselves with the case when the ``diffusion" coefficient of the set-valued SDE to be absolutely summable (see \S2 for details), mainly due to the path-regularity issue. We hope to address the well-posedness issue for more general (unbounded) set-valued SDE in our future publications.

The rest of this paper is organized as follows.  In \S2 we give the necessary preliminaries on set-valued analysis. 
In \S3 we introduce the $L_C$ space, a space of set-valued mappings with unbounded values and establish its basic properties. In \S4 we study the set-valued Lebesgue integrals on $L_C$-space and prove some important technical results including the continuity of the paths and some important estimates. 
In \S5  we prove the existence and uniqueness of strong solutions to a class of stochastic differential equations with   unbounded drift and absolutely summable  diffusion term. Finally, in \S6 we discuss the connections between the SDE and SDI in the unbounded set-valued case, and give an example of super-hedging problem with transaction costs in continuous time, which in fact motivated this work.

\section{Preliminaries}

In this section, we give a brief introduction to set-valued analysis and all the necessary
notations associated to it. Most of the materials can be found in the standard literature in set-valued analysis and stochastic analysis (see, e.g.,  \cite{BUC,MK,K}), but we shall try to give a self-contained presentation for ready references. 

%
Let $(\sX,d)$ be a metric (vector) space and denote $\sP(\sX):=2^{\sX}$ to be the family of all nonempty subsets of $\sX$,  and $\sC(\sX)\subset \sP(\sX)$ the family of all {\it closed} subsets of $\sX$. For a set $A \in \sP(\sX)$, we shall denote the closure of $A$ by $\overline{A}$ or $\cl A$, interchangeably. 
For $A,B \in \sP (\sX)$ and $\alpha \in \hR$, we define the   {\it Minkowski addition} and 
{\it scalar multiplication} by, respectively,  
 \bea
 \label{Minkovski}
 A+B:=\{a+b:a\in A , b\in B\}, \qq \alpha A:=\{\alpha a :a\in A\}.
 \eea
   It is important to note that in general, $A+(-1)A \neq \{\mathbf{0}\}$ unless $A$ is a singleton. Consequently, $\sP (\sX)$ is not a vector space under the Minkovski addition and scalar multiplication. 

Next, a set  $A \in \sP(\sX)$  is called   {\it convex} if for every $a_1,\cds,a_n\in A$ and  $\lambda_1,\cds, \lambda_n\in[0,1]$ such that $\sum_{i=1}^n\lambda_i=1$, it holds that $\sum_{i=1}^n\lambda_ia_i\in A$. We denote $\sG(\sX)$ to be the family of all {\it closed convex} sets. 
Let $A \in \sP(\sX)$ we define the {\it convex hull} of $A$, denoted by $\co A$, to be the smallest convex set containing $A$, and we denote  $\overline{\co} A:=\cl(\co A)\in \sG(\sX)$.
It is easy to check that the convex hull is exchangeable with the Minkovski addition and scalar multiplication. That is, if $A,B \in \sP (\sX) $ and   $\alpha \in\hR$, then it holds that
\bea
\label{coMinkov}
\co(A+B)=\co(A)+\co(B), \qq \co(\alpha A)=\alpha \ \co(A).
\eea

Finally, a  set $C\in \sP(\sX)$ is called a {\it cone} if $\alpha x \in C$ for any $\alpha>0 $ and $x\in C$. A  cone $C\in \sP(\sX) $ being  {\it closed} and {\it convex} are defined naturally. 
It is easy to check that  $C\in \sP(\sX)$  is a convex cone if and only if $C+C=C$ and $\alpha C=C$, for any $\alpha >0$. 

We should note that in this paper we shall focus on the space $\sC(\sX)$ or $\sG(\sX)$, which are non-compact in general. Therefore they are not closed under the Minkovski addition. For example, if $A, B\in \sC(\sX)$, then $A+B$ may not be a closed set, unless one of $A$ and $B$ is compact. A standard remedy for this is to introduce the extended notion of Minkovski addition, defined by
\bea
\label{Minkovexd}
A\oplus B:=\ol {A+B}, \qq A, B\in \sC(\sX).
\eea
Clearly, $A\oplus B=A+B$ if either $A$ or $B$ is compact. 

Let $A,B\in \sP(\sX)$. The {\it Hausdorff distance}  between $A$ and $B$   is defined by 
\bea
\label{Hausdorff}
h(A,B):=\max \{\overline{h}(A,B),\overline{h}(B,A) \}=\inf\{\e>0: A\subseteq B+\bar{B}_\e(0),~B\subseteq A+\bar{B}_\e(0)\}, 
\eea
 where $\bar{B}_r(0)\subset \hR^d$ is the closed ball centered at $0$ with radius $r>0$, and 
 $$\overline{h}(A,B):=\sup\{ d(a, B):a\in A\}:=\sup\{\inf\{d(a,b):b\in B\}:a\in A\}.$$ 
 It is clear by definition that $h(A,B)=0$ if and only if $\overline{A}=\overline{B}$, and it is well-known that
%
$(\sC (\sX),h)$ is a complete metric space if $\sX$ is. 
Furthermore, the following facts are used frequently in our future discussions: for $A,B,D,E\in \sP(\sX)$, then it holds that
\bea
\label{property}
\left\{\ba{lll}
h(\ol A,\ol B)\leq h(A,B); \qq 
h(\overline{\co} A,\overline{\co} B)\leq h(A,B); \ms\\ 
h(A\oplus B, D\oplus E)\le h(A+B,D+E)\leq h(A,D)+h(B,E)\ea\right.
\eea

\no {\bf Set-Valued Measurable Mappings and Decomposable Sets.}
 In what follows for simplicity we shall focus on the case $\sX=\hR^d$, the $d$-dimensional Euclidean space, although most of the results are valid for more general topological vector spaces. Let   $(\hX,\sM,\mu)$   be a finite measure space. If $\hX$ is a topological space, we take $\sM = \sB (\hX)$, the Borel $\sigma$-algebra on $\hX$. Consider a mapping $F:\hX\rightarrow \sP(\hR^d)$, which we shall refer to as the {\it set-valued mapping}.
A function $f:\hX\rightarrow \mathbb{R}^d$ such that $f(x)\in F(x) $, $x \in \hX$, is called a {\it selector} for $F$. A selector $f$ is called {\it measurable} if it is $\sM/\sB(\hR^d)$-measurable. 
We shall make use of the following definition of set-valued ``measurable" mapping.
\begin{defn}
A set-valued mapping $F:\hX\rightarrow \sP(\mathbb{R}^d)$ is said to be {measurable}  if for each $E\in \sC(\mathbb{R}^d)$, it holds that $\{x\in \hX:F(x)\cap E\neq \emptyset \}\in \sM$.
\qed
\end{defn}

The following selection/representation properties for measurable  mappings will be useful.
 
\begin{prop}[\cite{MK}]
\label{Fmeas}
(i) If $F:\hX\rightarrow \sC(\mathbb{R}^d)$ is measurable, then $F$ admits a measurable selector, i.e. there exists an $\mathcal{\sM}/\sB(\mathbb{R}^d)$-measurable selector $f:\hX\to \mathbb{R}^d$  of $F$; 

(ii) {\rm (Castaing Representation)} $F:\hX\rightarrow \sC(\mathbb{R}^d)$ is measurable if and only if there exists a sequence $\{f_n\}_{n=1}^\infty$ of  measurable selectors for $F$ such that $F(x)=\cl\{f_n(x):n\in \mathbb{N}\}$, $x\in \hX$.
\qed
\end{prop}

Let us denote  $\mathbb{L}^0(\hX,\mathbb{R}^d)=\mathbb{L}^0_{\sM}(\hX,\mathbb{R}^d)$ to be the set of all $\sM/\sB(\hR^d)$-measurable functions; and
$\mathcal{L}^0(\hX,\sC (\mathbb{R}^d))=\mathcal{L}^0_{\sM}(\hX,\sC (\mathbb{R}^d))$ the set of all measurable set-valued mappings. 
For $F\in \mathcal{L}^0(\hX,\sC (\mathbb{R}^d)),$ we consider the set of  measurable selectors
$$S(F):=\{f\in \mathbb{L}^0(\hX,\mathbb{R}^d):f(x)\in F(x)\text{~} \mu \text{-a.e.~} x\in \hX\}$$ 
Clearly, if $F,G\in \mathcal{L}^0(\hX,\sC (\mathbb{R}^d))$, then $F=G$, $\m$-a.e. if and only if $S(F)=S(G)$.
%

\begin{defn}
A set $V\subseteq \mathbb{L}^0(\hX,\mathbb{R}^d)$  is said to be  decomposable,    if    $\1_Af+\1_{A^c}g\in V$ whenever $f,g\in V$ and  $A\in\sM$.
\qed
\end{defn}

For a set $V\subseteq \mathbb{L}^0(\hX,\mathbb{R}^d)$, we define the  {\it decomposable  hull} of $V$, denoted by $\dec(V) $, to be the smallest decomposable set  in  $\mathbb{L}^0(\hX,\mathbb{R}^d)$ containing $V$. We shall often consider 
{\it closed decomposable  hull} of $V$:  $\overline{\dec}(V)=\cl[\dec(V)]$. It is not too hard to verify that the following algebraic properties hold for the decomposable hull: for $V,W\subseteq \mathbb{L}^0(\hX,\mathbb{R}^d)$ and $\alpha \in\hR$, 

\ss
1. $\dec(\alpha V)=\alpha \dec(V)$,

\ss
2. $\dec (V+W)=\dec(V)+\dec(W)$,

\ss
3. $\overline{\dec}[\co(V)]=\overline{\co}[\dec(V)]$. 

Let  $p\in[1,\infty)$ and  denote $\mathbb{L}^p(\hX,\mathbb{R}^d):=\{f\in \mathbb{L}^0(\hX,\mathbb{R}^d):\|f\|^p_p:=\int_X|f(x)|^p\mu(dx)<\infty\}$. For any $F\in \mathcal{L}^0(\hX,\sC (\mathbb{R}^d)),$ we define $S^p(F):=S(F)\cap \mathbb{L}^p(X,\mathbb{R}^d),$  and consider the set 
$$\sA^p(\hX,\sC (\mathbb{R}^d)):=\{F\in \mathcal{L}^0(\hX,\sC (\mathbb{R}^d)):S^p(F)\neq \emptyset\}.$$
We say $F$ is \textbf{$p$-integrable} if $F\in \sA^p(\hX,\sC (\mathbb{R}^d))$. 
 The following theorem (see \cite{MK}) will be important for our discussion. 
\begin{thm}
 \label{dec}
 (i) Let  $F\in  \sA^p(\hX,\sC (\mathbb{R}^d)).$ Then there exists a sequence $\{f_n\}_{n=1}^\infty\subseteq  S^p(F)$  such that $F(x)=\cl\{f_n(x):n\in \mathbb{N}\}$ for a.e. $x\in \hX$. Moreover, $S^p(F)=\overline{\dec}\{f_n:n\geq 1\}.$

%
\ms

(ii) Let $V\subset \mathbb{L}^p(\hX,\mathbb{R}^d)$ be a closed subset. Then there exists a measurable set-valued mapping $F\in \sA^p(\hX,\sC (\mathbb{R}^d))$ such that $V=S^p(F)$ if and only if $V$ is decomposable.
\qed
\end{thm}
Clearly, if   $F,G\in  \sA^p(\hX,\sC (\mathbb{R}^d))$ then $F=G~\mu$-a.e., if and only if $S^p(F)=S^p(G)$. Also, one can easily check that
\bea
\label{Sp}
 S^p\big(\overline{\co}F\big)=\overline{\co}S^p(F), \q
S^p\big(F\oplus G\big)= S^p(F)\oplus S^p(G).
\eea

%
%
%
%
Throughout the rest of this paper, we shall consider a given complete, filtered probability space $(\Omega,\mathcal{F},\mathbb{P},\mathbb{F}=\{\mathcal{F}_t\}_{t\in[0,T]}),$ on which is defined a standard $m$-dimensional Brownian motion $B=\{B_t\}_{t\in [0,T]}$, $T>0$. The notions of  set-valued random variables and stochastic processes, etc., can be defined in an analogues way as the usual ``single-valued" concepts in probability. For example, a {\it set-valued random variable}  $X:\Omega \rightarrow \sC (\mathbb{R}^d)$ is an $\mathcal{F}$-measurable set-valued mapping; and a {\it set-valued stochastic process} $\Phi=\{\Phi_t\}_{t\in[0,T]}$  is a family of set-valued random variables. A set-valued process is called {\it measurable} if it is $\sB([0,T])\otimes \mathcal{F}$-measurable; it is called {$\mathbb{F}$-adapted} if $\Phi_t$ is  $\mathcal{F}_t$-measurable  for each $t\in[0,T]$;  it is
called {\it $\mathbb{F}$-nonanticipative} if it is both measurable and adapted. 
Given a set-valued $\mathbb{F}$-non-anticipative process $\Phi,$ we denote  $S_{\mathbb{F}}(\Phi)$ to be the set of all $\mathbb{F}$-non-anticipative selectors of $\Phi$. For $p\ge 1$, we denote $S^p_{\mathbb{F}}(\Phi):=S_{\mathbb{F}}(\Phi)\cap  \mathbb{L}^p([0,T]\times \Omega,\mathbb{R}^d)$, and 
\bea
\label{AphF}
\sA^p_{\hF}([0,T]\times \O; \sC(\hR^d))=\{F\in \cL^0([0,T]\times \O;\sC(\hR^d)): S^p_\hF(F)\neq\emptyset\}.
\eea

We now turn out attention to the notions of set-valued Lebesgue and It\^{o} integrals. Let us assume that $\mathbb{F}=\mathbb{F}^B$ be the natural filtration generated by $B,$ augmented by all $\hP$-null sets of $\mathcal{F}$ so that it satisfies the usual hypothesis.

\textbf{Lebesgue Set-Valued Integral}:
Let $\Phi\in \sA^2_\hF([0,T]\times \O;\sC(\hR^d))$, and define
\bea
\label{J0T}
J _{0,T}( \Phi):=\big\{\int_0^T\phi_tdt:\phi\in S^2_{\mathbb{F}}(\Phi)\big\}\subseteq \mathbb{L}^2_{\mathcal{F}_T}(\Omega,\mathbb{R}^d).
\eea
 Since $\overline{\dec}_{\mathcal{F}_T}(J _{0,T} (\Phi))$ is closed and $\mathcal{F}_T$-decomposable,   by Theorem \ref{dec}, there exists a unique $\cF_T$-measurable set-valued mapping $I_T(\Phi)\in \sA^2_{\mathcal{F}_T}(\Omega,\sC (\mathbb{R}^d))$, such that $\overline{\dec}_{\mathcal{F}_T}(J _{0,T} (\Phi))=S^2_{\mathcal{F}_T}(I_T(\Phi))$.  We call $I_T(\Phi)$ the {\it set-valued Lebesgue integral} of $\Phi$, and denote it by $\int_0^T\Phi_tdt:=I_T(\Phi)$. Similarly, for $t\in [0,T]$, we  define $\int_0^t\Phi_sds:=\int_0^T\1_{(0,t]}\Phi_sds$, or equivalently, by  
 $\overline{\dec}_{\mathcal{F}_t}(J _{0,t} (\Phi))=S^2_{\mathcal{F}_t}(\int_0^t\Phi_sds)$. The following properties can be easily varified (cf., e.g., \cite{K, KisMic16}).

\ss 1. The set-valued integral $\int_0^t\Phi_sds$ is defined, almost surely, for each $t\in[0,T]$,  and is $\mathbb{F}$-adapted.

\ss
2. If $\Phi$ is convex-valued, so is $\int_0^t\Phi_sds$.

\ss
3.  If $\Phi$ is compact convex-valued, then the  process $\{\int_0^t\Phi_sds\}_{t\in[0,T]}$ is continuous with respect to the Hausdorff metric $h$.

We would like to point out here the continuity result of the Lebesgue integral $\{\int_0^t\Phi_sds\}_{t\in[0,T]}$ above  requires the compactness of $\Phi$. In fact, if $\Phi$ takes unbounded values,  the continuity is yet to be discussed in the literature, and will be one of the main focuses on this paper. 
%

\ss
\textbf{Aumann-It\^{o} Set-Valued Stochastic Integral}:
Let $\Psi\in \sA^2_{\hF}([0,T]\times\O;\sC (\mathbb{R}^{d\times m}))$. Define
\bea
\label{cJ0T}
\mathcal{J} _{0,T}(\Psi):=\Big\{\int_0^T\psi_tdB_t:\psi\in S^2_{\mathbb{F}}(\Psi)\Big\}\subseteq \mathbb{L}^2_{\mathcal{F}_T}(\Omega,\mathbb{R}^d).
\eea
 Again, 
by Theorem \ref{dec}, there exists a unique  set-valued mapping $\cI_T(\Psi)\in \sA^2_{\mathcal{F}_T}(\Omega,\sC (\mathbb{R}^d))$ such that $S^2_{\mathcal{F}_T}(\cI_T(\Psi))=\overline{\dec}_{\mathcal{F}_T}(\mathcal{J} _{0,T} (\Psi))$. We  denote $\cI_T(\Psi):=\int_0^T\Psi_tdB_t$, and call it the {\it Aumann-It\^{o} set-valued stochastic integral} of $\Psi$. Similarly,   we define $\int_0^t\Psi_sdB_s:=\cI_t(\Psi)$, $t\in [0,T]$, so that
 $S^2_{\mathcal{F}_t}(\cI_t(\Psi))=\overline{\dec}_{\mathcal{F}_t}(\mathcal{J} _{0,t} (\Psi) )$.
Again, the following facts can be easily verified (cf., e.g., \cite{K, KisMic16}):

\ss
1. The integral $\int_0^t\Psi_sdB_s$ is defined almost surely for each $t\in[0,T]$,  and it is $\mathbb{F}$-adapted. 

\ss
2. If $\Psi$ is convex-valued, so is $\int_0^t\Psi_sdB_s$.

However, the continuity of the process $\{ \int_0^t\Psi_sdB_s\}_{t\in[0,T]}$ as well as the joint-measurability is  more involved in the set-valued case. A well-known case is given by Kisielewicz-Mitchta \cite{KisMic16}, which we now describe. Let us  denote $\mathbb{L}^2_\hF([0,T]\times \Omega,\mathbb{R}^{d\times m})\subset\mathbb{L}^2([0,T]\times \Omega,\mathbb{R}^{d\times m})$ to be  the set of all non-anticipative processes, and let $G=\{g^n\}_{n\geq 1}\subseteq \mathbb{L}^2_\hF([0,T]\times \Omega,\mathbb{R}^{d\times m})$. Since for any $E\in\sC (\mathbb{R}^{d\times m})$,  it holds that
 $\{(t,\omega):G\cap E\neq\emptyset\}=\cup_{n=1}^\infty (g^n)^{-1}(E)$, 
 and each $g^n$ is  $\mathbb{F}$-non-anticipative, we see that $G$ is a $\mathbb{F}$-non-anticipative  set-valued process. Furthermore, 
 $G:=\{g^n\}_{n\geq 1}$
 is called {\it absolutely summable} if $\sum_{n=1}^\infty |g^n(t,\omega)|<\infty$ for $dt\times d\hP$-$(t,\omega)\in [0,T]\times \Omega$,  and it is called {\it square integrable} if $\mathbb{E}\int_0^T\sum_{n=1}^\infty |g^n(t,\cd)|^2 dt <\infty$. The following result gives a sufficient condition for the continuity of Aumann-It\^o set-valued stochastic integral:
  \begin{thm}[Kisielewicz-Mitchta \cite{KisMic16}]
\label{gcon} 
For every absolutely summable square integrable set-valued process
  $G:=\{g^n\}_{n\geq 1}\subseteq \mathbb{L}^2_\hF([0,T]\times \Omega,\mathbb{R}^{d\times m})$, the set-valued process   $\{\int_0^t G_sdB_s\}_{0\leq t \leq T}$ is almost surely continuous under Hausdorff metric $h$. 
   \end{thm}

We end this section by some technical results that will be  useful in our later sections.  Let $(\hX,\sM,\mu)$ be a finite measure space, and let $F\in  \sA^p_{\sM}(\hX,\sC(\hR^d))$ and $p\geq 1$.

\begin{prop}[{\cite[Theorem 2.3.4]{MK}}] 
\label{TechInf}
	 Let $\phi : \hX\times \hR^d\to [-\infty,\infty]$ be a jointly measurable function. Suppose that the integral $T_\phi (f):=\int_\hX \phi(x,f(x))\mu(dx)$ is well-defined for each $f\in S^p(F)$, and $T_\phi(f)< +\infty \mbox{~(resp. $>-\infty$)}$ for at least one $f\in S^p(F)$. Suppose further that either (i) $ \phi(x,\cd)$ is upper semicontinuous for every fixed $x\in \hX$ or  (ii) $(\hX,\sM,\mu)$ is a complete measure space and $\phi(x,.)$ is lower semicontinuous for every fixed $x\in \hX$. Then,
	\[
	\inf_{f\in S^p(F)} T_\phi(f)=\int_\hX \inf_{y\in F(x)} \phi(x,y)\mu(dx)\mbox{\rm ~(resp. $\dis\sup_{f\in S^p(F)} T_\phi(f)=\int_\hX \sup_{y\in F(x)}\phi(x,y)\mu(dx)$).} \qq\q
	\qed\]
	\end{prop}
%


Finally, we recall that a function  $f:\hX \times \mathbb{R}^d\rightarrow \mathbb{R}^d$ is said to be {\it Carath\`eodory}
if $f(x,\cd)$ is  continuous  for fixed $x\in \hX$, and $f(\cd,a)$ is measurable for  fixed $a\in  \mathbb{R}^d$. It can be shown that a  Carath\`eodory function $f:\hX \times \mathbb{R}^d\rightarrow \mathbb{R}^d$ must be $\sM\otimes \sB(\hR^d)$-measurable 
(see, e.g., {\cite[Theorem 6.1]{Him}}).  Accordingly, 
a set-valued mapping $F:\hX \times \mathbb{R}^d\rightarrow \sC(\mathbb{R}^d)$ is said to be a Carath\`eodory set-valued mapping if $F(x,\cd)$ is  continuous for  fixed $x\in \hX$,  and $F(\cd, a)$ is 
measurable for every fixed $a\in  \mathbb{R}^d$. Furthermore, we have the following result.
\begin{prop}[{\cite{SUP}}]
\label{FCAR}
Let $F:[0,T]\times \Omega\times \mathbb{R}^d\rightarrow \sC(\mathbb{R}^d)$ be a Carath\`eodory set-valued mapping, and let $F(t,\omega, A):=\overline{\cup_{a\in A}F(t,\omega,a)},$ for any $A\in \sC (\mathbb{R}^d).$  Then for every measurable $\mathbb{F}$-adapted set-valued mapping  $X:[0,T]\times \Omega\rightarrow  \sC (\mathbb{R}^d),$ the set-valued mapping $F\circ X:[0,T]\times \Omega\rightarrow  \sC(\mathbb{R}^d)$ such that  $(F\circ X)(t,\omega):=F(t,\omega, X(t,\omega))$ is $\sB([0,T])\otimes \mathcal{F}$ measurable. 
\qed
\end{prop}

\section{$L_C$-Space and Basic Properties}
\setcounter{equation}{0}

In this section we introduce the main subject of our framework, which we refer to as  the $L_C$-Space. As we observed  in the previous section, most of
the concepts of the standard set-valued analysis and stochastic analysis are defined on   $\sC(\hR^d)$ (or $\sG(\hR^d)$),  the space of all closed (convex) sets. However, many results in the set-valued analysis literature are valid only in $\sK(\hR^d)$, the space of all {\it compact}, convex sets. 
Our main task is  to remove the boundedness constraint in the analysis and establish a workable framework of set-valued analysis on unbounded sets. 

To begin with, we note that the first main issue  for allowing  unbounded sets is the choice of 
the reference element,  since  there is no single element in the space $\sC(\hR^d)$, to which the Hausdorff distance from
all elements is finite. Consequently, the commonly used reference point $\{\mathbf{0}\}$, as well as the 
``norm" $\|A\|=h(A, \{0\})$, is no longer valid. 
%
To  overcome this difficulty and proceed with our analysis, a reasonable remedy would be to consider a subspace of $\sC(\hR^d)$ for which there exists a fixed ``reference point" that has finite Hausdorff distance to all other elements. In light of the applications mentioned in the Introduction, in this paper we shall consider the  space $L_C$, which we now describe.

 Let  $ C \in \sG(\hR^d) $ be a fixed  convex cone in $\hR^d$, such that $\{\mathbf{0}\}\neq C\neq \hR^d$.  We consider the set: 
   \bea
\label{LC0}
L^\infty_C=L^\infty_C(\hR^d):=\{A\in\sC(\hR^d): A=\overline{A+C}\}.
\eea 
In what follows we do not distinguish $L^\infty_C$ and $L^\infty_C(\hR^d)$ when the context is clear. We endow the space $L^\infty_C$ with the usual Hausdorff metric $h$,  then it is not hard to check that  $L^\infty_C$ is ``closed" under countable union and intersections. Indeed,  if  $\{A_n\}_{n=1}^\infty\subseteq L^\infty_C$, then  $\overline{\cup_{n=1}^\infty A_n}\in L^\infty_C$ and $\cap_{n=1}^\infty A_n\in L^\infty_C$. We note that spaces similar to $L^\infty_C$ have appeared in the literature (e.g., \cite{setoptsurv}),  as a result of a preorder relation defined on linear spaces. In fact, we have the following result.
\begin{thm}
\label{LCspace}
$(L^\infty_C,h)$ is a complete metric space. 
\end{thm}

{\it Proof.}
Clearly, $(L^\infty_C,h)$ is a subspace of $(\sC(\hR^d), h)$. It suffices to show the completeness. To see this, let $\{A_n\}_{n=1}^\infty$ be a Cauchy sequence in $(L^\infty_C,h)$, whence a 
Cauchy sequence in $(\sC (\mathbb{R}^d),h)$. It then follows from \cite[Theorem 1.3.1]{MK} that  $\{A_n\}_{n=1}^\infty$ converges to $A:=\cap_{n=1}^\infty\overline{\cup_{m=n}^\infty A_m}$ in $(\sC (\mathbb{R}^d),h)$. Since each $A_m\in L^\infty_C$, so is $A$. 
\qed

We note that, however,  unlike the bounded case there is a main deficiency in the definition of $L^\infty_C$, that is, it is very possible for a  set $A\in L^\infty_C$ to have $h(A, C)=\infty$(!). Therefore, to facilitate our discussion, say, on SDE below, it is more desirable to consider the following subset of $L^\infty_C$:
\bea
\label{LC}
L_C=L_C(\hR^d):=\{A\in L^\infty_C(\hR^d):  h(A, C)<\infty\}.
\eea
Clearly, for all $A, B\in L_C$ we have $h(A, B)<\infty$,   thanks to the triangle inequality. However, we note that the space $L_C$ is no longer ``closed" under countable unions. The following simple example is thus worth noting:
\begin{eg}
\label{LCunion}
Let $A\in L^\infty_C(\hR^d)$ such that $h(A, C)=\infty$. For each $n\in\hN$, define $A_n:=(A\cap B_n(0))+C$, where $B_n(0)$ is the closed ball in $\hR^d$ centered at $0$ with radius $n$. Then, $A_n\in L_C$, as $h(A_n, C)\le n<\infty$, $n\in\hN$, but $A=\cup_{n=1}^\infty A_n \notin  L_C$.
\qed
\end{eg} 

We should note that, Example \ref{LCunion} notwithstanding, $L_C$ is still a complete metric space. 
\begin{prop}
\label{LCspace1}
The space $L_C$ is a closed subspace of $(L^\infty_C, h)$, hence a complete metric space. 
\end{prop}

{\it Proof.} We need only check that any limit point of $L_C$ belongs to $L_C$. To see this, let $(A_n)_{n\in\hN}\subset L_C$ and $\lim_{n\to\infty}h(A_n, A)= 0$.
Then $A\in L^\infty_C$, and for $N\in\hN$ large enough, we have $h(A_N, A)\le 1$. Thus we have $h(A, C)\le h(A, A_N)+h(A_N, C)<\infty$, since $A_N\in L_C$. To wit, $A\in L_C$. 
\qed

In the rest of the paper we shall focus only on the space $(L_C, h)$. Since the convex cone $C$ plays a special role in the space $L_C$, in what follows we refer to it as the {\it generating cone}.  Next, let us pay attention to a special type of elements in $L_C$. We say that 
 a set $B\in L_C$ has a {\it  compact component} $\tilde B\in \sK(\hR^d)$ if $B=\tilde B+C$. Clearly, if  in Example \ref{LCunion} $A\in L^\infty_C\cap \sG(\hR^d)$, then  the sequence $(A_n)$   
 all have compact components. An  important feature of such sets can be seen from the following extension of the so-called {\it Cancellation Law}\footnote{The ``Cancellation Law" states that for $A, B, C\in \sG(\hR^d)$, $A+B=C+B$ implies $A=C$, provided $B$ is compact (cf. \cite{PU}).}  to the case involving unbounded sets. 
\begin{prop}[Cancellation law]
\label{Cancel}
Let $A,B,D\in L^\infty_C \cap \sG(\hR^d) $ such that  $B$ has compact component.
Then $A\oplus B=D\oplus B$ implies that $A=D$.
\end{prop}

{\it Proof.} First assume that $ B= \tilde{B}+C$, where $\tilde{B}\in\sK(\hR^d)$. Since the Minkowski sum of a closed set and a compact set  is closed, we have $ A\oplus B = A\oplus (\tilde{B}+C)=\overline{A+(\tilde B+C)}=(A\oplus C)+\tilde{B}=A+\tilde{B}$. Similarly, we have $D\oplus B=D+\tilde{B}$. Therefore  $A\oplus B=D\oplus B$ implies that $A+\tilde{B}=D+\tilde{B}$. Since $\tilde{B}$ is compact, we conclude that $A=D$,  thanks to the usual cancellation law. 
\qed

We should note from Proposition \ref{Cancel} and Example \ref{LCunion} that the sets in $L_C$ that have compact components provides a useful subspace of $L_C$ for our discussion. In fact, such sets resemble the so-called   {\it Motzkin decomposable} sets in the literature  (see e.g. \cite{GGMT}). To see the relation between an element in $L_C$ that has compact component and a Motzkin decomposable set, we first recall  the so-called  {\it recession  cone} of a set $A\in \sP(\hR^d)$, defined by 
\bea
\label{recCone}
0^+A:=\{y\in \hR^d: x+\lambda y \in A,~ x\in A, ~\lambda \geq 0\}.
\eea
It is easy to see that $0\in 0^+A$, $A=A+0^+A$, and  $0^+A \neq \{0\}$ implies that  $A$ is unbounded. Moreover, if $A\in\sG(\hR^d)$ and $0^+A \neq \{0\}$,  then $0^+A$ is a non-trivial closed convex cone.  A set  $A\in \sG(\hR^d)$ with
$0^+A \neq \{0\}$ is  called {\it Motzkin decomposable} if $A =\tilde{A} +0^+A$ for some  $\tilde{A}\in \sK(\hR^d)$. 
We have the following result regarding the generating cone $C$ and the recession cone of $A\in L^\infty_C$.
\begin{prop}
\label{0+A}
For any $A\in L^\infty_C$, the following hold:

(i)  $0^+A \supseteq C$; 

(ii) If $A\in L_C$ and  $A$ has a compact component, then
 $A$ is  Motzkin-decomposable and $0^+A=C$; 
 
 (iii) If $A\in  L_C\cap \sG(\hR^d)$, then  $0^+A=C$.
\end{prop}

{\it Proof.}  (i) For any $c\in C$, $x\in A$, and $\lambda \geq 0$, it is readily seen that $x+\lambda c  \in A+ C\subseteq A$. 
Thus $C\subseteq 0^+A$. 
To see (ii), let $A\in L_C$, and let the compact component be $\tilde A$. 
Then we have $\tilde{A}+ C=A=A\oplus 0^+A=\tilde{A}+ (C \oplus 0^+A)$. Since $\tilde A$ is compact, by Proposition \ref{Cancel} we obtain $C= C\oplus 0^+A$. Now since $0\in C$, we have $ 0^+A=0^+A+\{0\}\subset 0^+A+C\subseteq C$, whence  $C= 0^+A$. 

%
It remains to check (iii).  First, by (i), $C\subseteq 0^+A$. On the other hand,  since $h(A, C)<\infty$, by definition of Hausdorff distance (\ref{Hausdorff}), we have
$ A\subseteq C+\bar{B}_\e(0)=:D_\e$, for some $\e>0$. Since $\bar{B}_\e(0)$ is compact, by (ii) $D_\e$ is Motzkin-decomposable and $0^+D_\e=C$. Now note that $A\subseteq D_\e$, by    definition (\ref{recCone}) we see that 
$C\subseteq 0^+A\subseteq 0^+D_\e=C$, proving (iii). 
\qed

We should note  that in general, $C\subsetneq 0^+A$ for $A\in L^\infty_C$. For example, let $C:=\{(x,y)\in \hR^2: x=0, y\geq 0\}$ and $A:=\{(x,y)\in \hR^2: x\geq 1, y\geq 1\}$. Then $A=A\oplus C$ but $0^+A=\hR^2_+$. 

To end this section let us consider the set-valued random variables taking values in the space $L_C$. Let us define 
\bea
\label{LC}
\left\{\ba{lll}
\mathbb{L}_C(\Omega,\mathcal{F};\hR^d):=\{F\in\mathcal{L}^0(\Omega,\sC(\mathbb{R}^d)) : F\in L_C, ~\hP\text{-a.s.}\}; \ms\\
\hL^2_C(\Omega,\cF;\hR^d):=\{F\in \hL_C(\Omega,\cF): \hE[h^2(F,C)]<\infty\}.
\ea\right.
\eea
We shall also often drop $\hR^d$ from the notations of the spaces above when the context is clear. Next, for any $F_1,F_2\in \hL^2_C(\Omega,\mathcal{F})$, define  $d(F_1,F_2):= (\mathbb{E}[h^2(F_1,F_2)] )^{\frac{1}{2}}$. 
 The following result is not surprising, except for some consideration due to the special structure of the $L_C$-space. We only give a sketch of proof for completeness.
\begin{thm}
\label{hLC}
$ (\mathbb{L}^2_C(\Omega,\mathcal{F}),d )$ is a complete metric space. 
\end{thm}

{\it Proof.} That $ (\mathbb{L}_C(\Omega,\mathcal{F}),d)$ is a  metric space is obvious, we shall argue only
%
%
%
the completeness. To this end, let   $(F_n)_{n=1}^\infty$ be a Cauchy sequence in $ (\mathbb{L}_C(\Omega,\mathcal{F}),d )$.  For any $\e >0$, by Chebyshev's inequality  we have $\hP\{h(F_n,F_m)>\epsilon \}\leq \frac{d^2(F_n,F_m) }{\epsilon ^2}$. Thus $(F_n)_{n=1}^\infty$
is Cauchy under Hausdorff distance $h$, in probability. 
We first claim that there is a subsequence $(F_{n_k})_{k=1}^\infty \subseteq (F_n)_{n=1}^\infty$ such that it converges $\hP$-a.s. in Hausdorff metric $h$ to some $F\in \hL_C(\Omega,\cF)$. 

Indeed, let $n_1:=1$ and for  $k>1$, 
$n_k:=\inf\{ n>n_{k-1}:\hP\{h(F_m,F_\ell)>2^{-k}\}<2^{-k}, ~m,\ell\geq n\}$.
 Let $E_k:= \{\omega\in \Omega : h(F_{n_{k+1}},F_{n_k})>2^{-k}\}$.  Then  $\sum_{k=1}^\infty \hP(E_k)<\sum_{k=1}^\infty 2^{-k}<\infty$. Denoting $\cN_l:=\cup_{k=l}^\infty E_k$ and $\cN:=\cap_{l=1}^\infty \cN_l$, we have $\hP(\cN)=0$, thanks to Borel-Cantelli lemma.  
 
 Now for any $l>1$, let $\omega\notin \cN_l$. Then for any $j>i>l$, it holds that
 $$h\big(F_{n_i}(\omega),F_{n_j}(\omega)\big)\leq \sum _{k=i}^{j-1}h\big(F_{n_{k+1}}(\omega),F_{n_k}(\omega)\big)\leq \sum _{k=i}^{j-1}2^{-k}\le 2^{-(i-1)}.$$
Thus $(F_{n_k}(\omega) )_{k=1}^\infty$ is a Cauchy sequence in $(L_C,h)$. 
Since  $(L_C,h)$ is complete,  by the proof of Theorem \ref{hLC}, for  $\omega\in \cN^c=\cup_{l=1}^\infty \cN_l^c$,  $(F_{n_k}(\omega) )_{k=1}^\infty$ converges to $\cap_{l=1}^\infty \overline{\cup_{k=l}^\infty F_{n_k}(\omega)}$ in $(L_C,h)$. 

Now let us define 
$F(\omega):=C\1_{\cN}(\o)+\cap_{l=1}^\infty \overline{\cup_{k=l}^\infty F_{n_k}(\omega)}\1_{\cN^c}(\o)$, $\o\in\O$, then $F(\omega)\in L_C$ for any $\omega \in \Omega$ and $(F_{n_k})_{k=1}^\infty$ converges to  $F$-$\hP$-a.s. We claim that $F$ is measurable and $(F_n)_{n=1}^\infty$ converges to $F$ in $(\hL_C(\Omega,\cF),d)$.  Indeed, for  $\emptyset \neq U\subseteq\hR^d$ open, we have
 \beaa
\{\omega\in\Omega:F(\omega)\cap U\neq \emptyset\}
=\{\omega\in\mathcal{N}:C\cap U\neq \emptyset\}\cup\{\omega\in \mathcal{N}^c:\cap_{l=1}^\infty \overline{\cup_{k=l}^\infty F_{n_k}(\omega)} \cap U\neq \emptyset\}.
\eeaa
Note that 
$\{\omega\in\mathcal{N}:C\cap U \neq \emptyset\}$ is a null set, 
hence measurable. On the other hand, 
$$\big\{\omega\in \mathcal{N}^c:\cap_{l=1}^\infty \overline{\cup_{k=l}^\infty F_{n_k}(\omega)}\cap U\neq \emptyset\big\}=\cap_{l=1}^\infty \cup_{k=l}^\infty\big\{\omega\in \mathcal{N}^c: F_{n_k}(\omega)\cap U\neq \emptyset\big\}$$
is also measurable. 
 Thus 
$\{\omega\in\Omega:F(\omega)\cap U \neq \emptyset\}$  is measurable, to wit,  $\o\mapsto F(\o)$ is measurable. 

Finally, we show that $(F_n)_{n=1}^\infty$  converges to $F$ in $(\hL_C(\Omega,\cF),d)$. Since $(F_n)_{n=1}^\infty$ is a Cauchy sequence in $(\hL_C(\Omega,\cF),d)$, then for any $\e  >0$ there is $N:=N(\e )\in\hN$ such that $\hE[h^2(F_n,F_m)]<\e $  for any $n,m\geq N$. Thus, for any $n>N$, we can apply Fatou's lemma to get
$$\hE[h^2(F_n,F)]=\hE[h^2(F_n,\lim_{{n_k}\rightarrow\infty}F_{n_k})]=\hE\big[\lim_{{n_k}\rightarrow\infty}h^2(F_n,F_{n_k})\big]\leq \liminf\limits_{{n_k}\rightarrow\infty}\hE[h^2(F_n,F_{n_k})]<\e. $$ 
This completes the proof. 
\qed

\section{Set-Valued Lebesgue Integral on $L_C$-Space }
\setcounter{equation}{0}

In this section, we show some results regarding set-valued Lebesgue integrals that will be useful for our discussion but not covered by the existing literature. We first recall the definition of set-valued Lebesgue integrals defined in \S2.
\begin{prop}
\label{DEC}
Let $F:[0,T]\times\Omega\rightarrow \sC(\hR^d)$ be a ``constant" set-valued mapping with $F\equiv A$, for some $A\in \sC(\hR^d)$. Then $F\in\sA^2_\hF([0,T]\times \Omega;\sC(\hR^d))$,  and for  $0\leq t_0<t\leq T $, it holds that $\int_{t_0}^tF_sds=(t-t_0)A$, $\hP$-a.s. 

In particular, if $A=C$ is a  convex cone, then  we have   $\int_{t_0}^tF_sds=C$, $\hP$-a.s. 
\end{prop}

{\it Proof.} Let $F(t, \o)\equiv A$, $(t, \o)\in[0,T]\times\Omega $ 
for some $A\in \sC(\hR^d)$. Then obviously $F$ is measurable and $\hF$-adapted.  
Let $\{a_n:n\geq 1\}:=A\cap\hQ^d$, then we have the simple Castaing representation: 
 $F(t,\omega)=\cl\{f^n_t(\omega):n\geq 1\}$, $(t,\omega)\in[0,T]\times \Omega$, where $f^n_t(\omega)\equiv a_n$. 
By (\ref{J0T}) we can easily see that
$$S^2_{\cF_t}\ ((t-t_0)A )=\overline{\dec}_{\cF_t} \{(t-t_0)a_n:n\geq 1\}= \overline{\dec}_{\cF_t} \Big\{\int_{t_0}^t f^n_sds :n\geq 1\Big\}=\overline{\dec}_{\cF_t}J_{0,T}(F), $$
 We conclude that $S^2_{\cF_t}(\int_{t_0}^tF_sds)=S^2_{\cF_t}((t-t_0)A)$ for any $0\leq t_0<t\leq T$. 
In particular, if $A=C$ is a convex cone, then we have   $(t-t_0) C=C$, proving the proposition.
\qed

Next, we give a standard result for Aumann-Lebesgue integral, extended to the unbounded integrands. Again, we only provide a sketch of proof for completeness.
\begin{prop}
\label{ADD}
  $F\in\sA^2_\hF([0,T]\times\Omega;\sC(\hR^d))$. Then for   $0\leq t_0<t\leq T $, it holds that
$$\int_0^tF_sds= \int_0^{t_0} F_sds\oplus\int_ {t_0}^t F_sds, \q \hP\as$$
\end{prop}
 {\it Proof.} For simplicity, let us denote  $I_0^t(F)=\int_0^tF_sds$, $I_0^{t_0}(F)=\int_0^{t_0}F_sds=\int_0^t\1_{[0, t_0 ]}F_sds$, and $I_{t_0}^t(F)=\int_{t_0}^tF_sds=\int_0^t\1_{[ t_0,t]}F_sds$, respectively. We shall argue that
\bea
\label{It0t}
S^2_{\cF_t} (I_0^t(F))=S^2_{\cF_t}  ( I_0^{t_0}(F)\oplus I _{ t_0}^t(F) ).
\eea 
Since by (\ref{J0T}),  $ S^2_{\cF_t}(I_0^t(F))= \overline{\dec}_{\cF_t}J_{0,t}(F):=\overline{\dec}_{\cF_t} \{\int_0^tf_sds\neg:f\in S^2_{\hF}(F) \}$, and by (\ref{Sp}), one has
\bea
\label{It0t1}
S^2_{\cF_t} ( I_0^{t_0}(F)\oplus I_{t_0}^t(F)  )= S^2_{\cF_t} (I_0^{t_0}(F) )\oplus S^2_{\cF_t} (I_{t_0}^t(F)).
\eea 
On the other hand, it is easy to check, regardless the boundedness of $F$, that
\beaa 
&& S^2_{\cF_t}(I_0^{t_0}(F))\oplus S^2_{\cF_t}(I_{t_0}^t(F)) = \overline{\dec}_{\cF_t}\{J_{0,t}(\1_{[0,t_0]}F) \}\oplus \overline{\dec}_{\cF_t}\{J_{0,t}(\1_{[t_0,t]}F) \}  \\
&=&  \dec_{\cF_t}\Big\{\int_0^tf_s^1ds:f^1\in S^2_{\hF}(\1_{[0,t_0]}F)\Big\}\oplus \dec_{\cF_t}\Big\{\int_0^tf_s^2ds:f^2\in S^2_{\hF}(\1_{[t_0,t]}F)\Big\} \\
&=& \overline{\dec}_{\cF_t}\Big\{\int_0^t(f_s^1+f_s^2)ds:f^1\in S^2_{\hF}(\1_{[0,t_0]}F), f^2\in S^2_{\hF}(\1_{[t_0,t]}F)\Big\}=\overline{\dec}_{\cF_t} \{J_{0,t}(F)\}.
\eeaa
Here, the last equality in the above can be easily verified using the definition of decomposibility. This, together with  (\ref{It0t1}), leads to (\ref{It0t}), proving the proposition.
\qed

 Next, we present an important  inequality regarding set-valued Lebesgue integrals. The  compact counterpart of this ineqaulity can be found in \cite{KisMic16}, we therefore only point out the special technical issues appearing in $L_C$-space. 
\begin{thm}
\label{FCON}
Suppose that $F,\tilde{F}\in\sA^2_{\hF}([0,T]\times\Omega;\sG(\hR^d))$ such that $\hE\big[\int_0^T h^2(F_s,\tilde{F}_s)ds \big] <\infty$.   Then, for any $0\leq t_0<t\leq T,$ it holds that
\bea
\label{h2ineq}
h^2\Big(\int_{t_0}^{t}F_sds,\int_{t_0}^{t}\tilde{F}_sds\Big)\leq (t-t_0)\int_{t_0}^th^2 (F_s,\tilde{F}_s )ds, \q \hP\as
\eea
\end{thm}

{\it Proof.} Fix $0\le t_0<t\le T$, and for notational simplicity, denote $I_{t_0}^t(\Phi):=\int_{t_0}^t\Phi_sds$, for $\Phi\in \sA^2_{\hF}([0,T]\times\Omega;\sG(\hR^d))$. Since both sides of (\ref{h2ineq}) is $\cF_t$ measurable, we need only to check that 
\bea
\label{h2ineq1}
\int_U h ^2 (I_{t_0}^{t}(F),I_{t_0}^{t}(\tilde{F}))d\hP \le (t-t_0)\int_U\int_{t_0}^th^2(F_s, \tilde{F}_s)dsd\hP, \q\forall U\in\cF_t.
\eea
To this end, let $U\in\cF_t$ be given, and consider the mapping $x\mapsto d^2(x, I_{t_0}^{t}(\tilde{F}))$. Since $d^2(\cd, E)$
is upper-semi continuous for  $E\in \hL_C(\O, \cF)$, it follows from Proposition \ref{TechInf} that
\bea
\label{LHS}
\int_U\overline{h}^2 (I_{t_0}^{t}(F),I_{t_0}^{t}(\tilde{F}))d\hP&=&\int_U\sup \{d ^2(x, I_{t_0}^{t}(\tilde{F})):x\in  I_{t_0}^{t}(F)\}d\hP \nonumber\\
&=&\sup\{\int_U d^2(u, I_{t_0}^t(\tilde{F}))d\hP:u\in S^2_{\cF_t}(I_{t_0}^{t}(F))\}.
\eea
By the similar argument, we can show that
\bea
\label{LHS1}
\int_Ud(u, I_{t_0}^{t}(\tilde{F}))d\hP=\inf\Big\{\int_U|u-v|^2d\hP:v\in S^2_{\cF_t}(I_{t_0}^{t}(\tilde{F}))\Big\}.
\eea 
Combining (\ref{LHS}) and (\ref{LHS1}) we obtain that
\bea
\label{LHS2}
\int_U\overline{h}^2(I_{t_0}^{t}(F),I_{t_0}^{t}(\tilde{F}))d\hP\neg\neg&\neg=\neg\neg&\neg\sup\Big\{\inf\Big\{\int_U|u-v|^2d\hP:v\in S^2_{\cF_t}(I_{t_0}^{t}(\tilde{F}))\Big\}:u\in S^2_{\cF_t}(I_{t_0}^{t}({F}))\Big\}\\
&=&\sup\Big\{\inf\Big\{\int_U|u-v|^2d\hP:v\in \overline{\dec}_{\mathcal{F}_t}J_{t_0,t}(\tilde{F})\Big\}:u\in  \overline{\dec}_{\mathcal{F}_t}J_{t_0,t}(F)\Big\}. \nonumber
\eea
Now note that by H\"older's inequality and Proposition \ref{TechInf} we also have 
\bea
\label{LHS3}
&&\sup\Big\{\inf\Big\{\int_U|u-v|^2d\hP:v\in J_{t_0,t} (\tilde{F}) \Big\}:u\in  J_{t_0,t} (F) \Big\}\nonumber\\
&=&\sup\Big\{\inf\Big\{\int_U\Big|\int_{t_0}^tf_sds-\int_{t_0}^t\tilde{f}_sds\Big|^2d\hP:\tilde{f}\in S^2(\tilde{F})\Big\}:f\in S^2(F)\Big\}\\
&\leq& 
(t-t_0)\sup\Big\{\inf\Big\{\int_{U\times [{t_0},t]}|f_s-\tilde{f}_s|^2d\hP ds:\tilde{f}\in S^2(\tilde{F})\Big\}:f\in S^2(F)\Big\}\nonumber\\
&=&(t-t_0)\int_U\Big[\int_{t_0}^t\overline{h}^2(F_s,\tilde{F}_s)ds\Big]d\hP . \nonumber
\eea
Now, in light of (\ref{LHS2}) we shall argue that (\ref{LHS3}) remains ture when we replace $J_{t_0,t} (F)$ and $J_{t_0,t} (\tilde{F})$ on the left hand side by their decomposable halls. To this end, first consider
$u=\sum_{i=1}^mu_i\1_{U_i}\in \dec_{\mathcal{F}_t} J_{t_0,t} (F) $,  where $u_i\in J_{t_0,t} (F) $ and $\{U_i\}_{i=1}^m$ is a  partion of $\O$. Then by (\ref{LHS3}) we see that
\bea
\label{LHS4}
&&\inf\Big\{\int_U|\sum_{i=1}^m u_i\1_{U_i}-v|^2d\hP:v\in \overline{\dec}_{\mathcal{F}_t }J_{t_0,t} (\tilde{F}) \Big\}\le \inf\Big\{\int_U|\sum_{i=1}^m u_i\1_{U_i}-v|^2d\hP:v\in J_{t_0,t} (\tilde{F}) \Big\}
\nonumber\\
&\le &\sum_{i=1}^m\sup\Big\{\inf\Big\{\int_{U\cap U_i}|u_i-v|^2d\hP:v\in J_{t_0,t} (\tilde{F}) \Big\}:u_i\in J_{t_0,t} (F) \Big\}\\
&\leq &\sum_{i=1}^m (t-t_0)\int_{U\cap U_i}\Big[\int_{t_0}^t \overline{h}^2(F_s,\tilde{F}_s) ds\Big]d\hP=  (t-t_0)\int_{U}\Big[ \int_{t_0}^t\overline{h}^2(F_s,\tilde{F}_s)ds\Big]d\hP . \nonumber
\eea
Since this is true for all $u\in \dec_{\mathcal{F}_t}J_{t_0,t} (F)$,  we see that (\ref{LHS4}) implies that (\ref{LHS3}) remains true when
$J_{t_0,t} (\tilde{F})$ is replaced by 
$\overline{\dec}_{\mathcal{F}_t }J_{t_0,t} (\tilde{F})$ and $J_{t_0,t} (F)$ is replaced by $\dec_{\mathcal{F}_t} J_{t_0,t} (F)$. To combine
(\ref{LHS2}) and (\ref{LHS3}), we let $u\in \overline{\dec}_{\mathcal{F}_t }J_{t_0,t} (F)$, and let  $(u_n)_{n=1}^\infty\subseteq \dec_{\mathcal{F}_t }J_{t_0,t} (F)$ such that $\lim_{n\rightarrow\infty}\hE|u-u_n|^2=0$.
Since for any $n\geq 1$, we can apply the extended (\ref{LHS3})  to get
\bea
\label{LHS5}
\inf\Big\{\int_U|u_n-v|^2d\hP :v\in \overline{\dec}_{\mathcal{F}_t }J_{t_0,t}( \tilde{F}) \Big\}\leq (t-t_0)\int_{t_0}^t\int_{U} \overline{h}^2(F_s,\tilde{F}_s)d\hP ds.
\eea
In particular,  (\ref{LHS5}) implies that for  $\e>0$ and $n\geq 1$,  we can find $v_n\in \overline{\dec}_{\mathcal{F}_t }J_{t_0,t}( \tilde{F} )$ such that 
\bea
\label{LHS6}
\int_U|u_n-v_n|^2d\hP\le 
 T\int_{ 0}^T\hE[\overline{h}^2(F_s,\tilde{F}_s)] ds+\e<\infty, 
\eea
thanks to the assumption of the theorem. Since $\{u_n\}_{n=1}^\infty $ is convergent in $ \hL^2([0,T]\times\O)$ by definition, whence bounded, this in turn implies that 
$\{v_n\1_{U}\}_{n=1}^\infty$ is bounded (in $\hL^2$) as well. We can then apply 
Banach-Saks-Mazur theorem to obtain a subsequence $(v_{n_k})_{k=1}^\infty$ of the form:
$\overline{v}_{n_k}:=\frac{1}{k}\sum_{j=1}^k v_{n_j}\1_U$, such that  $\{\overline{v}_{n_k}\}_{k=1}^\infty$ converges to, say, $\overline{v}\in\hL^2$. Since $\overline{\dec}_{\mathcal{F}_t }J_{t_0,t} (\tilde{F}) $ is convex, closed, and $\cF_t$-decomposable, we see that $\overline{v}_{n_k}\in \overline{\dec}_{\mathcal{F}_t }J_{t_0,t}( \tilde{F} )$, $k\geq 1 $, and hence
$\overline{v}\in \overline{\dec}_{\mathcal{F}_t }J_{t_0,t} (\tilde{F}) $. 
Moreover, since $\{u_n\}_{n=1}^\infty\subseteq \dec_{\mathcal{F}_t }J_{t_0,t} (F) $  converges to $u $, and
$\dec_{\mathcal{F}_t }J_{t_0,t} (F )$ is convex, one shows that $\overline{u}_{n_k}:=\frac{1}{k}\sum_{j=1}^ku_{n_j}\to u $ as well. Thus
by Jensen's inequality and (\ref{LHS5}) we have
\beaa
\int_U|\overline{u}_{n_k}-\overline{v}_{n_k}|^2d\hP &=&\int_U\Big|\sum_{j=1}^k\frac{1}{k}(u_{n_j}-v_{n_j})\Big|^2d\hP\le 
 \sum_{j=1}^k\frac{1}{k } \int_U|u_{n_j}-v_{n_j}|^2d\hP \\
&\le&(t-t_0)\int_{t_0}^t\int_{U} \overline{h}^2(F,\tilde{F})d\hP ds+\e.
\eeaa
Sending $k\to\infty$ we 
obtain that
$$\inf\Big\{\int_U|u-v|^2d\hP:v\in\overline{\dec}_{\mathcal{F}_t}J_{t_0,t}( \tilde{F}) \Big\} \le\int_U|u-\overline{v}|^2d\hP\leq (t-t_0)\int_{t_0}^t\int_{U} \overline{h}^2(F,\tilde{F})d\hP ds+\e.$$
Since $\e>0$ and $u\in \overline{\dec}_{\mathcal{F}_t}J_{t_0,t}(S^2(F))$ are arbitrary, and noting (\ref{LHS2}), we can first take a ``sup" in  the inequality above and then let $\e\to$ to conclude that (\ref{h2ineq1}) holds, proving the theorem.
\qed
\begin{rem}
A main technical point in this proof, compared to the compact case, is the application of Banach-Saks-Mazur theorem.  In fact, in the compact case, 
the set $\overline{\dec}_{\mathcal{F}_t}J_{t_0,t}( \tilde{F} ) $ is always bounded since 
 both $ \hE[\int_0^Th^2({F}_s,0)ds]$ and $\hE[\int_0^Th^2(\tilde{F}_s,0)ds] $ are finite. The assumption 
 $\hE\big[\int_0^T h^2(F_s,\tilde{F}_s)ds \big]<\infty$ is obviously much weaker. 
 \qed
 \end{rem}
 
Setting $\tilde{F}=C$ in Theorem \ref{FCON} and noting that $\int_{t_0}^t C=C$, we have the following corollary. 
 \begin{cor}
 \label{CCON}
Suppose that $F\in\sA_{\hF}^2([0, T]\times \O;\sG(\hR^d))$ such that $\hE\big[\int_0^Th(F_t, C)dt\big]<\infty$.  Then for any 
 $0\leq t_0<t\leq T$, it holds that
  $$h^2\Big(\int_{t_0}^{t}F_sds,C\Big)\leq (t-t_0)\int_{t_0}^th^2 (F_s,C )ds, \qq \hP\as
  $$
\end{cor}

To conclude this section, let us consider, 
similar to the spaces in (\ref{LC}), the following  spaces of $\hF$-nonanticipative set-valued mappings: 
\bea
\label{LCF}
\left\{\ba{lll}
\mathbb{L}_{C, \hF}([0,T]\times\Omega;\hR^d):=\{F\in\sA_\hF^2([0,T]\times\Omega;\sC(\hR^d)):F_t\in L_C, ~  \hP\text{-a.s.}, t\in[0,T] \}; \ms\\
\mathbb{L}^2_{C,\hF}([0,T]\times\Omega;\hR^d):=\{F\in \mathbb{L}_{C,\hF}([0,T]\times\Omega; \hR^d):\hE[\int_0^Th^2(F_s,C)ds]<\infty\}.  
\ea\right.
\eea
%
\begin{prop}
\label{Int}
 Let $F\in \mathbb{L}^2_{C,\hF}([0,T]\times\Omega;\hR^d)$. Then $\int_0^tF_sds\in \mathbb{L}^2_C(\Omega,\mathcal{F}_t)$, $t\in[0, T]$. Furthermore, the process $\{ \int_0^tF_sds\}_{0\leq t\leq T}$ has Hausdorff continuous paths, $\hP$-a.s.
\end{prop}
{\it Proof.}  Let $F\in  \mathbb{L}^2_{C,\hF}([0,T]\times\Omega)$. We first argue the adaptedness of the indefinite integral. Fix $t\in[0,T]$, by  definition of Lebesgue set-valued integral, it is clear that  $\int_0^tF_sds$ is an $\cF_t$-measurable set-valued random variable. Furthermore, since $C$ is a convex cone containing $\bf0$, we have  $\int_0^tCds=Ct=C$. This, together with the fact that
 $F_t\in  L_C$, $\hP$-a.s. $t\in[0,T]$, yields
 $$\int_0^tF_sds\oplus C=\int_0^t F_sds\oplus \int_0^tCds=\int_0^t (F_s\oplus C)ds=\int_0^tF_sds.
$$
 Moreover, by defintion (\ref{LCF}) and Corollary \ref{CCON} we have $\hE[h^2(\int_0^tF_sds ,C)]\leq t \hE[\int_0^th^2(F_s,C)ds]<\infty$.
 Thus $\int_0^tF_sds\in \mathbb{L}^2_C(\Omega,\mathcal{F}_t)$, $t\in[0,T]$.

To see the continuity, let $0\leq t_0<t\leq T.$ By   Theorem \ref{DEC}, Propositions \ref{ADD}, \ref{Int}, Corollary \ref{CCON},   and properties of Hausdorff metric, we have, $\hP$-almost surely,
\beaa
h^2\Big(\int_0^tF_sds,\int_0^{t_0}F_sds\Big)&=&h^2\Big( \int_0^{t_0}F_sds\oplus \int_{t_0}^tF_sds, \int_0^{t_0}F_sds\oplus C \Big)\\
&\leq& h^2\Big(\int_0^{t_0}F_sds+ \int_{t_0}^tF_sds,\int_0^{t_0}F_sds + C\Big)\\
&\le &h^2 \Big(\int_{t_0}^tF_sds,C\Big)\leq (t-t_0) \int_0^T h^2(F_s,C)ds.
\eeaa
Thus $h^2 (\int_0^tF_sds,\int_0^{t_0}F_sds ) \rightarrow 0$, as $|t-t_0|\rightarrow 0$, proving the result. 
\qed

\section{Set-Valued SDEs with Unbounded Coefficients}
\setcounter{equation}{0}

In this section, we consider stochastic differential equations of the form:
\bea
\label{SDE}
X_t =  \xi\oplus \int_0^t F(s,X_s)ds\oplus \int_0^t \mathcal{G}\circ XdB_s, \qq t\in[0,T], ~\hP\as,
\eea
 where $\xi\in \mathbb{L}^2_C(\Omega, \mathcal{F}_0)$ and $F, \mathcal{G}$ are set-valued mappings taking valued in $\sC(\hR^d)$ (hence possibly unbounded), which we now describe.  
%
\begin{defn}
\label{F}
A mapping $F:[0,T]\times \Omega\times \mathbb{R}^d\rightarrow L_C $ is called a non-anticipating {\it Carath\'eodory} set-valued random field if it enjoys the following properties: 

(i)  $F\in \cL^0( [0,T]\times\O\times \hR^d;L_C)$;

(ii) For fixed $a\in\hR^d$, $F(\cd, \cd, a)\in \hL_{C,\hF}([0,T]\times\O;\hR^d)$; and 

(iii) For $\hP$-a.s. $\o\in\O$, $F(\cd, \o, \cd)$ is a Carath\`eodory  set-valued function.
%

Furthermore, for a given non-anticipating Carath\'eodory set-valued random field $F$, we define, with a slight abuse of notations, its set-to-set version $F:[0,T]\times \Omega\times L_C  \rightarrow L_C $ by 
\bea
\label{coefficientF}
F(t,\omega, A):=\overline{\co}(\cup_{a\in A}F(t,\omega,a), \qq A\in L_C, ~ (t, \o)\in[0,T]\times\O.
\eea
 \end{defn}
Clearly, if the coefficient  $F$  in SDE (\ref{SDE}) is a non-anticipating Carath\'eodory set-valued random field of the form (\ref{coefficientF}), then it follows from Proposition \ref{FCAR} that for any $X\in L^0_{C,\hF}([0,T]\times \Omega; \hR^d)$, $F(\cd, X_\cd)\in \cL^0_\hF([0,T]\times\O;\hR^d)$. Furthermore, we make the following assumption.
\begin{assum}
\label{assumpF}
$F:[0,T]\times\O\times L_C\to L_C$ is a non-anticipative Carath\'eodory random field such that 
for some constants $\beta>0$,  the following conditions hold:


(1) $h^2(F(t,\cd,A),C)\leq \beta (1+ h^2(A,C)) $, $A\in L_C$, $\hP$-a.s.;

(2) $h^2(F(t,\cd,A),F(t,\cd,\tilde{A}))\leq \beta h^2(A,\tilde{A})$, $A,\tilde{A}\in L_C$, $\hP$-a.s.
\qed
\end{assum}

Next, we specify the coefficient $\cG\circ X$ in (\ref{SDE}). We note that in SDE (\ref{SDE}) the stochastic integral is defined in Aumann-It\^o sense, and we shall try to define the process $\cG\circ X$ so that Theorem \ref{gcon} can be applied. 
We therefore begin with the following definition.
\begin{defn}
\label{coefficientg}
A mapping  $g:[0,T]\times \Omega\times L_C\rightarrow \mathbb{R}^{d\times m}$ is called a non-anticipating {\it Carath\'eodory} random field if it enjoys the following properties: 

(i) $g\in \hL^0( [0,T]\times\O\times L_C; \hR^{d\times m})$;

(ii) for fixed $A\in L_C$, $g(\cd, \cd,A)\in \hL^0_{\hF}([0,T]\times\O;\hR^{d\times m})$; and 

(iii) for fixed $(t, \o)$, the mapping $A\mapsto g(t,  \o, A)$,  is Hausdorff continuous. 
%
 \end{defn}
Now let us consider the class of sequences $\{g^n\}$ of non-anticipating Carath\'eodory random fields that further satisfy the following uniform Lipschitz condition.
\begin{assum}
\label{assumpG}
There exist  $\alpha_n >0$, $n=1,2, \cds$,  such that $\sum_{n=1}^\infty \a^2_n<\infty$, and 

(1)  $|g^n(t,\cd,C)| \leq \alpha_n$, $\hP$-a.s.;

(2) $|g^n(t,\cd,A)-g^n(t,\omega,\tilde{A})| \leq \a_n h(A,\tilde{A})$,  $A,\tilde{A}\in L_C$, $\hP$-a.s. 
%
\qed
\end{assum}

We remark that Assumption \ref{assumpG} implies the following growth condition.
\bea
\label{gngrowth}
|g^n(t, \cd, A)|\le \a_n(1+h(A, C)), \q A\in L_C.
\eea
Furthermore, if $G:=\{g^n\}$ is a sequence of non-anitcipating Carath\'eodory random fields satisfying Assumption \ref{assumpG}, and 
$X\in \hL^2_{C, \hF}([0,T]\times \Omega; \hR^{d\times m})$, then in light of Theorem \ref{FCAR}, for each $n$ the mapping 
$g^n\circ X \in\hL^0_\hF([0,T]\times \Omega; \mathbb{R}^{d\times m})$, where $(g^n \circ X)(t, \o):=g^n(t, \o, X(t,\omega))$, $(t, \o)\in[0,T]\times\O$. Furthermore, since $X\in \hL^2_{C, \hF}([0,T]\times \Omega; \hR^{d\times m})$, we have $\hE[\int_0^T h^2(X_t, C)dt]<\infty$, 
thus
\bea
\label{Gest}
\hE\Big[\int_0^T \sum_{n=1}^\infty |(g^n\circ X)(t, \cd)|^2dt\Big]&\le& \hE\Big[\int_0^T \sum_{n=1}^\infty \a^2_n [1+h^2(X(t, \cd), C)]dt\Big]\\
&=&\Big(\sum_{n=1}^\infty \a^2_n \Big) \hE\Big[\int_0^T [1+h^2(X(t, \cd), C)]dt\Big]<\infty. \nonumber
\eea 
In what follows, for $G=\{g^n\}$ and $X\in \hL^2_{C, \hF}([0,T]\times\O;\hR^{d})$ 
as above we denote 
$\cG\circ X :=\overline{\co}\{g^n\circ X:n\geq 1\}=\ol{\co}(G\circ X)$, then  the Aumann-It\^o integral 
$\int_0^t \cG\circ X dB_s=\ol{\co}[\int_0^t G\circ X dB_s]$ is well-defined for each $t\in[0,T]$, and it follows from Theorem \ref{gcon} that the process of indefinite integrals $\{\int_0^t\cG\circ X dB_s\}_{t\in[0,T]}$ has Hausdorff continuous paths. 
%

Our main result of the well-posedness for set-valued SDE (\ref{SDE}) is the following theorem.
\begin{thm} 
\label{existuniq}
Assume that the Assumptions \ref{assumpF} and \ref{assumpG} hold for the coefficients $F$ and $\cG$, respectively. Then, for each $\sG(\hR^d)$-valued random variable $\xi\in\mathbb{L}^2_C(\Omega,\mathcal{F}_0 )$, there exists a unique convex, (Hausdorff) continuous  process $X=(X_t)_{0\leq t\leq T}\in \mathbb{L}^2_{C, \hF}([0,T]\times\Omega;\hR^d)$, such that (\ref{SDE}) holds $\hP$-a.s.
\end{thm}

{\it Proof.} We follow the standard Picard iteration: let 
$Y^{(0)}\equiv \xi$, and for $k\ge 0$, we define 
\bea
\label{Yk}
Y_t^{(k+1)}:=\xi\oplus \int_0^t F(s,Y_s^{(k)})ds\oplus \int_0^t \cG \circ Y^{(k)}dB_s, \qq t\in[0,T].
\eea
We first claim that for each  $k\ge 0$,  $Y^{(k)}\in \hL^2_{C, \hF}([0,T]\times \O; \hR^d)$  and has Hausdorff continuous paths. Indeed, since $\xi\in \hL^2_C(\O, \cF_0; \hR^d)$, the claim is trivial for $k=0$. Now assume that the claim is true for $Y^{(k)}$. Then, by  using Assumptions \ref{assumpF}, \ref{assumpG}, along with the estimates similar to (\ref{Gest}) and an induction argument, one can easily check that $F(\cd, Y^{(k)}_\cd) \in\mathbb{L}^2_{C, \hF}([0,T]\times \Omega;\hR^d)$  and   $\cG\circ Y^{(k)}\in \hL^2([0,T]\times\O;\hR^{d\times m})$. Then by Proposition \ref{Int} and 
Theorem \ref{gcon}, we further deduce that all the indefinite  integrals $\int_0^tF(s,Y^{(k)})ds$ and $\int_0^t\mathcal{G}\circ Y^{(k)}dB_s$, $t\in [0,T]$, are well-defined $\hF$-adapted set-valued processes with Hausdorff continuous paths. Thus so is  $Y^{(k+1)} $.
%

To show that $Y^{(k+1)}\in \hL^2_{C, \hF}([0,T]\times \O; \hR^d)$, it remains to check $h(Y^{(k+1}_t, C)<\infty$, $\hP$-a.s. To see this, 
let $t\in[0,T]$ be fixed, then noting (\ref{property}) and the fact that $C$ is a cone, we have
\bea
\label{est1}
&&\hE[h^2(Y^{(k+1)}_t,C)]=\hE\Big[h^2\Big(\xi\oplus\int_0^t F(s,Y^{(k)}_s)ds\oplus\int_0^t \mathcal{G}\circ Y^{(k)}dB_s\},C\Big)\Big]\nonumber\\
&\leq& \hE\Big[h^2\Big(\xi+\int_0^t F(s,Y^{(k)}_s)ds+\int_0^t \mathcal{G}\circ \xi dB_s,C+C+\{0\}\Big)\Big]\\
&\leq& 2 \hE\Big[h^2\Big(\xi+
\int_0^t F(s,Y^{(k)}_s)ds,C+C\Big)\Big]+2 \hE\Big[h^2\Big(\int_0^t \mathcal{G}\circ Y^{(k)}dB_s,\{0\}\Big)\Big]\nonumber\\
&\leq& 4 \hE [h^2 (\xi,C ) ]+4\hE\Big[h^2\Big(\int_0^t F(s,Y^{(k)}_s)ds,C\Big)\Big]+2 \hE\Big[h^2\Big(\int_0^t \mathcal{G}\circ Y^{(k)}dB_s,\{0\}\Big)\Big]. \nonumber
\eea
Note that, by Corollary \ref{CCON} and the inductional assumption, we have 
\bea
\label{est2}
\hE\Big[h^2\Big(\int_0^t F(s,Y^{(k)}_s)ds,C\Big)\Big]\leq t \hE\Big[\int_0^t h^2(F(s,Y^{(k)}_s),C)ds\Big]<\infty. 
\eea
Moreover,  by  \cite[Theorem 5.4.2]{MK}  and Assumption \ref{assumpG}, we have 
\bea
\label{est3}
\mathbb{E}\Big[h^2\Big(\int_0^t \mathcal{G}\circ Y^{(k)} dB_s,\{\mathbf{0}\}\Big)\Big]\leq \mathbb{E}\Big[\int_0^t\sum_{n=1}^\infty |g^n\circ Y^{(k)}_s|^2 ds\Big]<\infty.
\eea
Combining (\ref{est1})-(\ref{est3}) we obtain $\hE[h^2(Y^{(k+1)}_t,C)]<\infty$, $t\in[0,T] $, whence the claim. 

Next, let us estimate $\hE[h^2(Y^{(k+1)}_t,Y^{(k)}_t)]$, for $k\ge 0$. Again, recall (\ref{property}), we have 
\bea
\label{est4}
&&\mathbb{E} [h^2 (Y_t^{(k+1)},Y_t^{(k)})]\\
&\leq& \mathbb{E}\Big[h^2\Big( \xi+\int_0^t F(s,Y_s^{(k)})ds+\int_0^t \mathcal{G}\circ Y^{(k)}dB_s, \xi+\int_0^t F(s,Y_s^{(k-1)})ds+\int_0^t \mathcal{G}\circ Y^{(k-1)}dB_s\Big)\Big]\nonumber\\
&\leq& 2 \mathbb{E}\Big[h^2\Big(\int_0^t F(s,Y_s^{(k)})ds,\int_0^t F(s,Y_s^{(k-1)})ds\Big)\Big]+2 \mathbb{E}\Big[h^2\Big(\int_0^t \mathcal{G}\circ Y^{(k)}dB_s,\int_0^t \mathcal{G}\circ Y^{(k-1)}dB_s\Big)\Big].\nonumber
\eea
Now, applying Theorem \ref{FCON} and using the Assumption \ref{assumpF}, we have 
\bea
\label{est5}
&&\mathbb{E}\Big[h^2\Big(\int_0^t F(s,Y_s^{(k)})ds,\int_0^t F(s,Y_s^{(k-1)})ds\Big)\Big] \leq  T\mathbb{E}\Big[\int_0^t h^2 (F(s,Y_s^{(k)}), F(s,Y_s^{(k-1)})ds\Big)ds\Big]\nonumber\\
&\leq &T\beta\mathbb{E}\Big[\int_0^t  h^2(Y_s^{(k)},Y_s^{(k-1)})ds\Big]=T\beta \int_0^t \mathbb{E} [ h^2(Y_s^{(k)},Y_s^{(k-1)})]ds.
\eea
Furthermore, by  \cite[Theorem 5.4.2]{MK}) and Assumption \ref{assumpG} we obtain 
\bea
\label{est6}
&&\mathbb{E}\Big[h^2\Big(\int_0^t \mathcal{G}\circ Y^{(k)} dB_s,\int_0^t \mathcal{G}\circ Y^{(k-1)}dB_s\Big)\Big]
\leq  \mathbb{E}\Big[\sum_{n=1}^\infty\Big| \int_0^t [g^n\circ Y^{(k)} -g^n\circ Y^{(k-1)}]dB_s\Big|^2 \Big]\nonumber\\
&\le&  \mathbb{E}\Big[\sum_{n=1}^\infty \int_0^t\a_nh^2(Y^{(k)} ,Y^{(k-1)})ds \Big]
= \Big(\sum_{n=1}^\infty \a_n\Big) \int_0^t\mathbb{E}[ h^2(Y^{(k)} ,Y^{(k-1)})]ds. 
\eea
Combining (\ref{est4})-(\ref{est6}) we obtain 
\bea
\label{est7}
\mathbb{E}\big[h^2\big(Y_t^{(k+1)},Y_t^{(k)}\big)\big]&\leq & 2\Big[T\beta+\sum_{n=1}^\infty  \a_n\Big]\int_0^t\hE\big[h^2\big(Y_s^{k},Y_s^{(k-1)}\big)\big]ds.
\eea
Since  $Y{(0)}\equiv \xi$, by Assumptions \ref{assumpF} and \ref{assumpG} we can easily check that
\bea
\label{est8}
\mathbb{E}\big[h^2\big(Y_t^{(1)},Y_t^{(0)}\big)\big] \le    2\Big[\beta T+\sum_{n=1}^\infty  \alpha_n \Big]  (1+\mathbb{E}[h^2(\xi,C)] )t, 
\qq t\in[0,T]. 
\eea
Repeatedly applying (\ref{est7}) and noting (\ref{est8}), an induction argument leads to that
\bea
\label{est9}
\mathbb{E}\big[h^2\big(Y_t^{(k+1)},Y_t^{(k)}\big)\big] \le  M^{k+1}
(1+\mathbb{E}[h^2(\xi,C)]) \frac{t^{k+1}}{(k+1)!}, 
\eea
where $M:= 2\big[T\beta+\sum_{n=1}^\infty  \a_n\big]$.
Now consider the complete metric space  $(\hL_C(\O, \cF_t), d)$ in Theorem \ref{hLC}. In terms of the metric $d$ we see that (\ref{est9}) implies that,  for $m>n $, 
\beaa
d(Y_t^{(m)},Y_t^{(n)})&\leq &\sum_{k=n}^{m-1}d(Y_t^{(k+1)},Y_t^{(k)})= \sum_{k=n}^{m-1}\big(\mathbb{E}\big[h^2(Y_t^{(k+1)},Y_t^{(k)})\big]\big)^{\frac{1}{2}}\\
&\leq &(1+\mathbb{E}[h^2(\xi,C)])^{\frac12} \sum_{k=n}^{m-1}\Big(\frac{M^{k+1}t^{k+1}}{(k+1)!}\Big)^{\frac{1}{2}}
\rightarrow 0, \q\text{~as~} m, n\rightarrow \infty.
\eeaa
Thus $\{Y_t^{(k)}\}_{k\in\mathbb{N}}$ is a Cauchy sequence in   $(\mathbb{L}_C(\Omega,\mathcal{F}_t),d)$. Thus  there exists $Y_t\in \mathbb{L}_C(\Omega,\mathcal{F}_t)$ such that $\lim_{k\rightarrow \infty}d(Y_t,Y_t^k)=0$. We shall argue that the limit $Y $ has Hausdorff continuous paths. More precisely, we claim that  $\{Y^{(k)}\}_{k\in\mathbb{N}}$ convergences  to $Y$ uniformly on $[0,T]$. Indeed, denoting $\D^{(k)}g^n(t, \cd):=g^n(t,Y_t^{(k)} )-g^n(t,Y_t^{(k-1)})$, $\D^{(k)}(h\circ F)(t, \cd):=h(F(t,Y_t^{(k)} ),F(t,Y_t^{(k-1)})$, $k\in\hN$, and applying Theorem \ref{FCON}  we can easily show that, for any $k\ge 0$, $\hP$-almost surely, 
\bea
\label{est10}
&&\sup_{0\leq t\leq T}h^2(Y_t^{(k+1)},Y_t^{(k)}) \leq   2 \sup_{0\leq t\leq T} t \int_0^t [\D^{(k)} (h\circ F)]^2(s,\cd )ds
+2 \sup_{0\leq t\leq T}\sum_{n=1}^\infty \Big|\int_0^t |\D^{(k)}g^n(s,\cd)]dB_s\Big|^2\nonumber\\
&&\qq\qq \leq 2T\beta  \int_0^T h^2 (Y_s^{(k)},Y_s^{(k-1)} )ds+2 \sup_{0\leq t\leq T}\sum_{n=1}^\infty \Big|\int_0^t |\D^{(k)}g^n(s,\cd)|dB_s\Big|^2.
\eea
Here the second inequality above is due to Assumption \ref{assumpF}. Thus, for each $k\in\hN$ we have
\bea
\label{est11}
&&\hP\Big(\sup_{0\leq t\leq T}h^2(Y_t^{(k+1},Y_t^k)>\frac{1}{2^k}\Big)\le  \hP\Big(2T\beta  \int_0^T h^2 (Y_s^{(k)},Y_s^{(k-1)})ds >\frac{1}{2^k}\Big)\nonumber\\
&&\qq\qq\qq+ \hP\Big(2\sup_{0\leq t\leq T}\sum_{n=1}^\infty \Big|\int_0^t [\D^{(k)}g^n(s,\cd)]dB_s\Big|^2>\frac{1}{2^k}\Big).
\eea
Now, by Markov's inequality we have
\bea
\label{est12}
\hP\Big(2T\beta  \int_0^T h^2 (Y_s^{(k)},Y_s^{(k-1)})ds >\frac{1}{2^k}\Big) \leq  2^{k+1}T\beta  \mathbb{E}\Big[ \int_0^T h^2(Y_s^{(k)},Y_s^{(k-1)} )ds \Big].
\eea
Furthermore, by Burkholder-Davis-Gundy's inequality, for some $K>0$, it holds that
\bea
\label{est13}
&&\hP\Big(2\sup_{0\leq t\leq T}\sum_{n=1}^\infty \Big|\int_0^t [\D^{(k)}g^n(s,\cd)]dB_s\Big|^2>\frac{1}{2^{k}}\Big) 
%
%
%
%
%
%
\le 2^{k+1}\mathbb{E}\Big[ \sum_{n=1}^\infty \sup_{0\leq t\leq T}\Big|\int_0^t[\D^{(k)}g^n(s,\cd)]dB_s\Big|^2\Big]\nonumber\\
&\le & 2^{k+1}K\sum_{n=1}^\infty \mathbb{E}\Big[ \int_0^T |\D^{(k)}g^n(s,\cd)|^2ds\Big] \le 2^{k+1}K\Big(\sum_{n=1}^\infty  \a_n\Big) \mathbb{E}\Big[ \int_0^Th^2(Y_s^{(k)},Y_s^{(k-1)})ds\Big].
\eea
Combining (\ref{est11})-(\ref{est13}) and noting (\ref{est9}) we deduce that
\bea
\label{est14}
&&\hP\Big(\sup_{0\leq t\leq T}h^2(Y_t^{k+1},Y_t^k)>\frac{1}{2^k}\Big)\nonumber\\
&\leq& 2^{k+1}T\beta  \mathbb{E}\Big[ \int_0^T h^2\Big(Y_s^k,Y_s^{k-1}\Big)ds \Big]+ 2^{k+1}K\Big(\sum_{n=1}^\infty  \a_n\Big) \mathbb{E}\Big[ \int_0^Th^2(Y_s^{(k)},Y_s^{(k-1)})ds\Big]\nonumber\\
&=& 2^{k+1}\Big[T\beta+K\sum_{n=1}^\infty  \a_n\Big]\int_0^T\mathbb{E} [ h^2 (Y_s^{(k )},Y_s^{(k-1)}) ]ds\\
&\le & 2^{k+1}\Big[T\beta+K\sum_{n=1}^\infty  \a_n\Big](1+\hE[h^2(\xi, C)])\frac{M^kT^{k+1}}{(k+1)!}. \nonumber
\eea
Therefore, we conclude that $\sum_{k=1}^\infty \hP(\sup_{0\leq t\leq T}h^2(Y_t^{(k+1)},Y_t^{(k)})>\frac{1}{2^k} )<\infty$. A standard argument 
using the Borel-Cantelli lemma then yields that the sequence
$\{Y_t^{(k)}\}_{k\in\mathbb{N}}$ converges  uniformly on $[0,T]$ to the  $Y\in\hL^2_{C, \hF}([0,T]\times \O)$, $\hP$-a.s.,  and thus $Y$ has (Hausdorff) continuous paths.
%
%

Finally, the similar argument shows that, for fixed $t\in[0,T]$, 
\beaa
&&\lim_{k\rightarrow\infty}\mathbb{E}\Big[h^2\Big(\int_0^tF(s,Y_s^{(k)})ds,\int_0^tF(s,Y_s)ds\Big)\Big]
+\lim_{k\rightarrow\infty}\mathbb{E}\Big[h^2\Big(\int_0^t\mathcal{G}\circ Y^{(k)}dB_s,\int_0^t\mathcal{G}\circ YdB_s\Big)\Big]\\
&\leq&\Big[T\beta+ \sum_{n=1}^\infty \a_n\Big] \lim_{k\rightarrow\infty} \mathbb{E}\Big[ \int_0^t h^2 (Y_s^{(k)},Y_s)ds\Big]=0.
\eeaa
That is, $Y $ satisfies the SDE (\ref{SDE}). The argument of the uniqueness can be carried out in a similar manner, we leave it to the interested reader. 
\qed
%
%
%

To end this section, we give the following {\it stability} results of the  solution to the  SDE (\ref{SDE}). 
\begin{prop}
\label{stable}
Assume that Assumptions \ref{assumpF} and \ref{assumpG} are in force, and let $X$ be the solution to SDE (\ref{SDE}). Then for any $0\leq s<t\leq T$, it holds that
\bea
\label{stabest}
\left\{\ba{lll}
\hE[\sup_{0\leq t\leq T}h^2(X_t,C)]\leq K_T  \Big(1+\mathbb{E}[h^2(\xi,C)]\Big); \ms\\
\hE[h^2(X_t,X_s)]\leq   K_T (1+ \mathbb{E} [h^2(\xi,C)])(t-s).
\ea\right.
\eea
where $K_T>0$ is a generic constant depending only on $T$, $\beta$ in Assumption \ref{assumpF}, and $\{ \alpha_n\}$ in Assumption 
\ref{assumpG}.
\end{prop}

{\it Proof.} The proof is fairly standard by now. We only give a sketch. First, in what follows let us denote $K_T$ to be a 
generic constant depending only on $T$, $\beta$, and $\{\alpha_n\}$, which is allowed to vary from line to line. 
Then, similar to (\ref{est1}) we have, for $t\in[0,T]$, 
\bea
\label{est15}
&&\mathbb{E}\Big[\sup_{s\leq t}h^2(X_s,C)\Big] =
\mathbb{E}\Big[\sup_{s\leq t}h^2\Big(\Big\{  \xi\oplus\int_0^s F(\tau,X_\tau)d\tau\oplus \int_0^s \mathcal{G}\circ XdB_\tau\Big\},C\Big)\Big]\\
&\leq& 4\Big\{ \mathbb{E}\Big[h^2(\xi,C)\Big]+ \mathbb{E}\Big[\sup_{s\leq t}h^2\Big(\int_0^sF(\tau,X_\tau)d\tau,C\Big)\Big]+ \mathbb{E}\Big[\sup_{s\leq t}h^2\Big(\int_0^s(\cG\circ X)_\tau dB_\tau,\{\mathbf{0}\}\Big)\Big]\Big\}. \nonumber
\eea

By Theorem \ref{FCON} and Assumption \ref{assumpF} we have
\bea
\label{est16}
&&\mathbb{E}\Big[\sup_{s\leq t}h^2\Big(\int_0^sF(\tau,X_\tau)d\tau,C\Big) \leq  \mathbb{E}\Big[\sup_{s\leq t}s\int_0^s h^2\Big(F(\tau,X_\tau),C\Big)d\tau\Big]\\
&\leq& K_T\mathbb{E}\Big[\int_0^t (1+ h^2(X_\tau,C))d\tau \Big] 
\le  K_T \Big(1+\mathbb{E}\Big[\int_0^t \sup_{s\leq \tau} h^2(X_s,C)d\tau \Big]\Big). \nonumber
\eea
Similarly, by Assumption \ref{assumpG} and Burkholder-Davis-Gundy inequality we have
\bea
\label{est17}
&&\hE\Big[\sup_{s\leq t}h^2\Big(\int_0^s(\cG\circ X)  dB_\tau,\{\mathbf{0}\}\Big)\Big]\leq \mathbb{E}\Big[\sup_{s\leq t}\sum_{n=1}^\infty \Big|\int_0^s (g^n\circ X)_\tau dB_\tau\Big|^2\Big]\\
& \leq&  \sum_{n=1}^\infty  \mathbb{E}\Big[\int_0^t |(g^n\circ X)_\tau |^2d\tau\Big] 
\le K_T\Big(1+\mathbb{E}\Big[\int_0^t \sup_{s\leq \tau}h^2(X_s,C)d\tau\Big]\Big).  \nonumber
\eea
Combining (\ref{est15})-(\ref{est17}) we obtain
\beaa
\mathbb{E}\Big[\sup_{s\leq t}h^2(X_s,C)\Big] 
\le K_T\Big(1+\int_0^t \mathbb{E}\Big[\sup_{s\leq \tau}h^2(X_s,C)\Big]d\tau\Big), \qq t\in[0,T].
\eeaa
The first inequality in (\ref{stabest}) then follows from the Gronwall inequality. To see the second inequality of (\ref{stabest}) we note that, 
for $0\le s<t\le T$, 
%
\bea
\label{est18}
&&\hE[h^2(X_t,X_s)] \nonumber\\
&=&\hE\Big[h^2\Big(\int_0^tF(\tau,X_\tau)d\tau\oplus \int_0^t(\cG\circ X )dB_\tau,\int_0^sF(\tau,X_\tau)d\tau\oplus\int_0^s(\cG\circ X)dB_\tau \Big)\Big]\\
&\leq &2 \hE\Big[ h^2\Big(\int_0^tF(\tau,X_\tau)d\tau,\int_0^sF(\tau,X_\tau)d\tau\Big) \Big]+2 \hE\Big[ h^2\Big( \int_0^t(\cG\circ X)dB_\tau,\int_0^s(\cG\circ X)dB_\tau\Big) \Big]. \nonumber
\eea
Now applying Proposition \ref{Int} and noting the first inequality of (\ref{stabest}) we have
\bea
\label{est19}
&&\hE\Big[ h^2\Big(\int_0^t F(\tau,X_\tau)d\tau,\int_0^s F(\tau,X_\tau)d\tau \Big)\Big] \nonumber\\
& \leq &\hE\Big[ h^2\Big(\int_0^s F(\tau,X_\tau)d\tau+\int_s^tF(\tau,X_\tau)d\tau,\int_0^s F(\tau,X_\tau)d\tau+C\Big)\Big] \\
&\leq &\hE\Big[ h^2 \Big(\int_s^tF(\tau,X_\tau)d\tau,C\Big)\Big] 
\leq  T\hE\Big[  \int_s^th^2(F(\tau,X_\tau),C)d\tau\Big]
 \leq  (t-s)K_T \big(1+ \mathbb{E} [h^2(\xi,C) ]\big). \nonumber
\eea
Similarly, we have
\bea
\label{est20}
&& \hE\Big[ h^2\Big( \int_0^t(\cG\circ X)dB_\tau,\int_0^s(\cG\circ X)dB_\tau\Big) \Big] 
  \leq   \sum_{n=1}^\infty \hE\Big[\int_s^t|g^n(\tau,X_\tau)|^2 d\tau\Big]\nonumber\\
&\leq &K_T \hE\Big[\int_s^t\big(1+h^2(X_\tau,C)) d\tau\Big] 
\leq (t-s)  K_T \Big(1+\hE\Big[\sup _{0\le \tau\leq T}h^2(X_\tau,C)\Big]\Big)\\
&\leq &(t-s)   K_T  \big( 1+\mathbb{E} [h^2(\xi,C) ]\big).\nonumber
\eea
Combining (\ref{est18})-(\ref{est20}) we derive the second inequality of (\ref{stabest}). The proof is now complete. 
\qed

\section{Connections to SDIs and Applications}
\setcounter{equation}{0}

The theoretical framework we have established in the previous sections can be used to study the {\it Stochastic Differential Inclusions} (SDI) with unbounded coefficients which, to the best of our knowledge, has not been systematically investigated. In fact, the interplay between the solutions to SDIs and those of  SVSDEs, as well as some applications of SDIs in the so-called continuous time {\it Super-hedging} theory was one of the main motivations of this study. In this section we shall give a brief discussion on these issues, and we refer to our forthcoming work \cite{ACM} for a more detailed study of the aspects in finance. 

Let us begin with the description of an SDI associated with the SVSDE (\ref{SDE}). Let $F:[0,T]\times\O\times \hR^d\rightarrow L_C$ and $\{g^n:n\geq 1\}$ be a family of functions $g^n:[0,T]\times \O\times \hR^d\rightarrow \hR^{d\times m}$ satisfying the following {\it Standing Assumptions}:


\begin{assum}
\label{assumpSDI}
The functions $F$ and $\{g^n\}$ are  non-anticipative Carath\'eodory random fields such that: 

(i) there exists $\beta>0$, such that for   fixed $t\in [0,T]$, $x, y\in \hR^d$,  and $ A\in L_C$, it holds that
\bea
\label{LipF}
\left\{\ba{lll}
\sup_{x\in A}h^2(F(t,\o, x),C)\le \beta (1+h^2(A,C)),   \ms\\
h^2(F(t,\o, x),F(t,\o, y))\leq \beta |x-y|,\ea\right. \qq \hP\ae \o\in\O.
\eea

(ii) there exist constants $\a_n>0$ such that $\sum_{n=1}^\infty \a^2_n<\infty$, and for   fixed $t\in [0,T]$ and $x, y\in \hR^d$, it holds that
\bea
\label{Lipgn}
|g^n(t, \o, x)|\le \a_n, \q |g^n(t,\o, x)-g^n(t,\cd, y)|\leq \a_n|x-y|,  \qq\hP\ae \o\in\O.
\eea
\end{assum}
Assuming that $F$ and $\{g^n\}$ are Carath\'eodory (set-valued) random fields satisfying Assumption \ref{assumpSDI}, and for a given process $\bx\in\hL^2_\hF([0,T]\times \O;\hR^{d\times m})$, we  define $(\cG\circ \bx)_t :=\ol{\co}\{(g^n\circ \bx)_t:
n\ge 1\}$. Then, by Theorem \ref{gcon}, the Aumann-It\^{o} indefinite integral $\{\int_0^t(\cG\circ \bx)dB\}$ is well-defined and has (Hausdorff) continuous paths. We can now consider the the SDI in the following sense.
\begin{defn}
\label{defnSDI}
A process $\bx=(\bx_t)_{0\leq t \leq T}\in \hL^2_\hF([0,T]\times\O;\hR^d)$ is said to be a solution to a Stochastic Differential Inclusion if

\ss
(i) $\bx$ has continuous paths; and 

\ss
(ii) the following relations holds: 
\bea
\label{SDI}
\bx_t-\bx_0\in \int_0^t F(s,\bx_s)ds\oplus \int_0^t (\cG\circ \bx) dB_s, \qq 0< t \leq T, \q \hP\as
\eea
\end{defn}

The SDI (\ref{SDI}) is closely related to the  SDE studied in the previous section. For example, let $F$ and $\{g^n\}$ be the Carath\'eodory random fields satisfying Assupmtion \ref{assumpSDI}, and let us define 
$$\ol{F}(t,\o, A):=\overline{\co}[\cup_{x\in A}F(t,\o, x)], \q \ol{g}^n(t,\o, A):=\sup_{x\in A} g^n(t,\o, x), \q (t, \o)\in[0,T]\times \O, ~A\in L_C.
$$
Here $\sup_{x\in A} g^n(\cd, \cd, x):=\big(\sup_{x\in A}g^n_{ij}(\cd, \cd, x)\big)_{i,j=1}^{d, m}$. 
We claim that 
the mappings $\ol{F}$ and $\ol{G}:=\{\ol{g}^n\}$ are Carath\'eodory set-valued random fields satisfying Assumptions \ref{assumpF} and \ref{assumpG}, respectively. Indeed, note that by Assumptions \ref{assumpSDI}-(2) we have
\bea
\label{olF1}h^2(\ol{F}(t,\o, A) ,C)\leq h^2(\cup_{x\in A}F(t,\o, x),C)\leq \sup_{x\in A} h^2(F(t,\o,x),C)\leq \beta (1+ h^2(A,C))<\infty.
\eea
Furthermore, since $F(t,\o, x)\in L_C$, we have $F(t,\o, x)=F(t,\o, x)\oplus C$. Thus, 
\bea
\label{olF2}
\ol{F}(t,\o, A)=\overline{\co}[\cup_{x\in A}(F(t,\o, x)\oplus C)]
=\overline{\co}  (\cup_{x\in A} F(t,\o, x)) \oplus C=\ol{F}(t,\o, A)\oplus C.
\eea
Here in the above we used the fact that $\ol\co(C)=C$. Clearly, (\ref{olF1}) and (\ref{olF2}) imply that $\ol{F}$ is a non-anticipative Carath\'eodory $L_C$-valued random field defined on $[0,T]\times\O\times L_C$, and that Assumption \ref{assumpF}-(1) holds.  
%
Moreover, Assumption \ref{assumpF}-(2)  follows from (\ref{LipF}) and \cite[Lemma (4.1,1)]{KisMic16}. To see Assumption \ref{assumpG}, we note that (\ref{Lipgn}) implies that
$|\overline{g}^n(t,\cd,A)| \leq \sup _{x\in A}|g^n(t,\cd,x)|\leq \a_n$, $\hP$-a.s., whence Assumption \ref{assumpG}-(1). To vefify Assumption \ref{assumpG}-(2), we first note that all Euclidean norms on $\hR^{d\times m}$ are equivalent, so we shall use the ``maximum" norm on $\hR^{d\times m}$,  that is,  $|M|:=\max_{i,j}|M_{i,j}|$, $M\in\hR^{d\times m}$. Now let $A, \tilde{A}\in L_C$, then, for any $x\in A$, $y\in\tilde{A}$, by (\ref{Lipgn}) we have
 \beaa
\label{bargn1}
g^n_{ij}(t, \o, x)\le g^n_{ij}(t, \o, y) +\a_n|x-y|\le \bar{g}^n_{ij}(t, \o, \tilde{A})+\a_n|x-y|,\q 1\le i\le d, ~1\le j\le m.
\eeaa
Since $x\in A$,  $y\in\tilde{A}$, and $i,j$ are arbitrary, it is readily seen that the inequality above yields
\bea
\label{bargn2}
\bar{g}^n_{ij}(t, \o, A)=\sup_{x\in A}g^n_{ij}(t, \o, x)\le \bar{g}^n_{ij}(t, \o, \tilde{A})+\a_n\bar{h}(A, \tilde{A}).
\eea
Switching the position of $A$ and $\tilde{A}$ in (\ref{bargn2}) and recalling the definition of the maximum norm we see that 
$\{\bar{g}^n\}$ satisfies Assumption \ref{assumpG}-(2). 
%
%

We can then consider the following set-valued SDE:
\bea
\label{SDESDI}
X_t= \xi\oplus \int_0^t\ol{F}(s,X_s)ds\oplus \int_0^t(\ol{\cG}\circ X)dB_s, \qq t\in [0,T], 
\eea
which is well-posed, thanks to Theorem \ref{existuniq}. 
%
%
\begin{rem}
\label{remark6.0}
{\rm The connection between the solutions to SDI (\ref{SDI}) and the selectors of the solution to SDE (\ref{SDESDI}) is an interesting issue. In the compact coefficient case, it is known that there exists at least one solution to the SDI that is also a selector of the solution to SDE (cf. \cite{KisMic16}). It would be interesting to see if an analogue for SDI (\ref{SDI}) and SDE (\ref{SDESDI}) remains true, but the arguments is actually quite subtle due to the unboundedness (whence non-compactness) of the coefficients. We prefer to pursue it in our future publications. 
\qed}
\end{rem}
%

To end this section we shall 
give an example in finance that actually motivated this work. More precisely, we extend some well-known concepts regarding {\it super-hedging problems with transaction costs} (cf. e.g., \cite{cvi,Kcurrency, Sch}) to a continuous-time framework, and show how this would lead to an SDI whose drift coefficients will contain the so-called {\it solvency cone}. whence $L_C$-valued, and thus unbounded. 
We begin by a market model proposed in \cite{cvi}, which we now describe.
 
Consider a financial market consisting of one risky asset $S$ and one riskless asset  $B$ whose dynamics are given by 
\bea
\label{BSSDE}
\left\{\ba{lll}
dS_t=S_t[b_tdt+\sigma_t dW_t], \qq &S(0)=p\in(0,\infty);\\
dB_t=r_tB_t dt, &B_0=1.
\ea\right.
\eea
A {\it trading strategy} is a pair $(L,M)$ of $\hF$-adapted, nondecreasing processes with $L_0=M_0=0$, where $L_t$ (resp. $M_t$) represent  the total amount of funds transferred from $B$ to $S$ (resp. from $S$ to $B$). 
We shall consider the case where there exist (proportional) transition costs in both directions, and denote $0<\lambda<1$ to be  the proportion of the cost   from  $B$ to $S$, and $0<\mu<1$ be that   from $S$ to $B$. Assuming that all the transaction costs are charged/depsited to the bank (riskless asset) account, then    given  initial holdings $(x,y)$   and a trading strategy $(L,M)$, we denote {\it portfolio process}, representing the amount of funds in the riskless and risky assets, by $X =X^{x,L,M}$ and $Y =Y^{y,L,M} $, respectively. Then $(X, Y)$ should have the dynamics: for $0\leq t \leq T$, 
\bea
\label{CVIXY}
\left\{\ba{lll}
dX_t= r_tX_tdt-(1+\lambda)dL_t+(1-\mu)dM_t, \qq &X_0=x;\\
dY_t=Y_t [b_t dt+\sigma_t dW_t] +dL_t-dM_t , &Y_0=y.
\ea\right.
\eea

Let us now describe the problem in terms of the {\it Solvency Cone} proposed in \cite{Kcurrency, KandS} and \cite{Sch}. We begin by 
considering two assets, still denoted by $B$ and $S$, respectively, but without requiring $B$ being riskless. In order that our argument can be extended to more assets, let us denote  $\pi^{12}_t$ to be the number of units of $B$ that can be exchanged for one share of $S$ at time $t$ (including its transaction cost). That is, we have, $\hP$-a.s.,
\bea
\label{pi12}
\pi^{12}_t B_t  =(1+\lambda) S_t  \qq\mbox{or} \qq \pi^{12}_t =(1+\lambda)  \frac{S_t }{B_t }, \qq t\in[0,T].
\eea
Similarly, we let  $\pi^{21}_t$ denote the number of shares of $S$ that can be exchanged  for  one unit of $B$  at time $t$ (including its transaction cost), so that the following identities hold: 
\bea
\label{pi21}
\pi^{21}_t S_t  =\frac{1}{1-\mu} B_t  \qq\mbox{or}\qq \pi^{21}_t =\frac{1}{1-\mu} \frac{B_t }{S_t }, \qq t\in[0,T], \q \hP\as
\eea
Now let us define $\pi^{11}_t=\pi^{22}_t :\equiv1$ and consider the matrix-valued process $\Pi_t:=(\pi^{ij}_t)_{i,j=1}^2$, $t\in[0, T]$. Clearly, the components $\pi^{ij}_t$ satisfy the following properties:
%
%
\bea
\label{bidask}
\pi^{ij}_t >0,  i,j =1, 2; \qq   \pi^{11}_t =\pi^{22}_t =1;\qq \pi^{12}_t \pi^{21}_t >1, \qq t\in[0,T], \q \hP\as
\eea
We shall follow \cite{Kcurrency, Sch} and call any matrix $\Pi$ satisfying (\ref{bidask}) the {\it bid-ask matrix}, and the process $\Pi=(\Pi_t)_{0\leq t\leq T}$ taking values in bid-ask matrices the {\it  bid-ask process}. Next we defined the so-called {\it Solvency Cone} associated with the bid-ask process $\Pi=(\Pi_t)_{0\leq t\leq T}$. Recall that for a given set of vectors $\{\xi^1, \cds, \xi^n\}\subset \hR^d$, 
a cone  generated by $\{\xi^i\}_{i=1}^n$ is define by 
$ K:=\text{cone~}\{\xi^1, \cds, \xi^n\}=\{\sum_{i=1}^n \a_i\xi^i: \a_i\ge 0\}.
$
Now for a given bid-ask process $\Pi_t=\{\pi^{ij}_t\}_{i,j=1}^2$, $t\in[0,T]$, we define
	\bea
	\label{KPi}
	K(\Pi_t ):=\text{cone}\{e^1,e^2,\pi^{12}_t e^1-e^2,\pi^{21}_t  e^2-e^1\} =\text{cone}\{(1,0),(0,1),(\pi^{12}_t,-1),(-1,\pi^{21}_t ) \}.
\eea

Note that by (\ref{pi12}) and (\ref{pi21}),  $\pi^{12}_t=(1+\l)\frac{S_t}{B_t}>0$, and $\pi^{21}=\frac1{1-\m}\frac{B_t}{S_t}>0$, we see that
$(0, 1), (1,0)\in \text{cone}\{(\pi^{12}_t,-1),(-1,\pi^{21}_t )\}$, therefore, using the conic properties we can easily deduce that
\bea
\label{KPi1}
K(\Pi_t )\neg\neg&\neg\neg=\neg\neg&\neg\neg \text{cone} \Big\{\neg\Big((1+\lambda)  \frac{S_t}{B_t} ,\neg-1\Big),\Big(\neg-\neg1,\frac{1}{1-\mu} \frac{B_t}{S_t}\Big)\Big\}
\neg=\neg\text{cone} \Big\{\neg\Big(\neg(1+\lambda)  \frac{1}{B_t} ,-\frac{1}{S_t}\Big),\Big(\neg-\neg(1-\mu)\frac{1}{B_t}   ,\frac{1}{S_t}\Big)\neg\Big\}\nonumber\\
&=&\Big\{\Big((1+\lambda)  \frac{\a}{B_t}-(1-\mu)\frac{\beta }{B_t}, \frac{\beta-\a}{S_t}\Big): ~~ \a, \beta \ge 0 \Big\}.
\eea

Now let us change the unit of the portfolio process $(X, Y)$ in (\ref{CVIXY}) to ``number of shares" and denote
$h^1_t=X_t/B_t$, $h^2_t=Y_t/S_t$, $t\in[0,T]$. We shall call $h_t=(h^1_t, h^2_t)$, $t\in[0,T]$, the {\it unit portfolio process}, and we have the following ``inclusion" result.
\begin{prop}
Assume that the portfolio process $(X, Y)$ satisfies (\ref{CVIXY}), and the price process $(B, S)$ satisfies (\ref{BSSDE}). Assume further that there exists a   rate process $\th=(\th^L_t, \th^M_t)_{t\in[0, T]}$, such that $\th^L_t\ge 0, \th^M_t\ge 0$, $\hP$-a.s., and that  the trading strategy $(L, M)$ satisfies $ L_t=\int_0^t\th^L_rdr$ and $M_t=\int_0^t\th^M_sds$, $t\in[0,T]$. Then, the corresponding unit portfolio process $h=(h^1, h^2)$ satisfies the following relation $\hP$-almost surely:
\bea
\label{hcone}
h_t=(h^1_t, h^2_t)\in (x,y/p)+\int_0^t [ -K(\Pi_s)]ds, \qq t\in[0, T].
\eea
\end{prop}

{\it Proof.} First, by It\^{o}'s formula and some simple calculation, we have 
\bea
\label{dh}
\left\{\ba{lll}
dh^1_t=d (
X_tB^{-1}_t)
= B^{-1}_t[-(1+\lambda)dL_t+(1-\mu)dM_t\Big], \qq & h^1_0=x; \\
dh^2_t= d ( Y_t S^{-1}_t)
=  S^{-1}_t [-dM_t+dL_t], & h^2_0=y/p.
\ea\right.
\eea
Since $d(L_t, M_t)=(\th^L_t, \th^M_t)dt $ by assumption, we can write (\ref{dh}) as  
\bea
\label{hdynamic}
h_t&=&(h^1_t, h^2_t)=(x,y/p)-\Big((1+\lambda)\int_0^t\frac{\th^L_s}{B_s}ds -(1-\mu)\int_0^t\frac{\th^M_s}{B_s}ds,\int_0^t \frac{\th^M_s-\th^L_s}{S_s}ds\Big)\nonumber\\
&=&(x,y/p)-\int_0^t\Big((1+\lambda)\frac{\th^L_s}{B_s} -(1-\mu) \frac{\th^M_s}{B_s}, \frac{\th^M_s-\th^L_s}{S_s}\Big)ds, \q t\in[0,T].
\eea
Since $\th^L_t\ge 0$, $\th^M_t\ge0$, $\hP$-a.s., $t\in[0,T]$, by (\ref{KPi1}) we see that 
$$\Big((1+\lambda)\frac{\th^L_t}{B_t} -(1-\mu) \frac{\th^M_t}{B_t}, \frac{\th^M_t-\th^L_t}{S_t}\Big)\in K(\Pi_t), ~~\hP\mbox{\rm -a.s.,~~}t\in[0,T].
$$
Thus (\ref{hcone}) follows from (\ref{hdynamic}).
\qed

 To derive the SDI of type (\ref{SDI})  let us now consider a deterministic cone:
\bea
\label{K}
K:=\Big\{ \Big(1+\lambda)\alpha -(1-\mu)\beta, \beta -\alpha\Big): \alpha, \beta \geq 0\Big\}\subseteq \hR^2.
\eea
Then it is easy to see from (\ref{KPi1}) that $(x,y)\in K$ if and only if $(\frac{x}{B},\frac{y}{S})\in K(\Pi)$. Consequently, since $\th^L_t\ge 0$, 
$\th^M_t\ge 0$, we have 
 $$\big( (1+\lambda)\th^L_t-(1-\mu)\th^M_t,\th^M_t-\th^L_t \big)\in K,~~\hP\mbox{\rm -a.s.,~~} t\in[0,T].$$
We can now rewrite (\ref{CVIXY}) as 
\beaa
\label{XYdynamic}
(X_t,Y_t)& =&(x,y)-\Big( (1+\lambda)L_t-(1-\mu)M_t,M_t-L_t \Big)+\Big(\int_0^tX_s r_sds, \int_0^tY_s[b_sds+\sigma_sdW_s]\Big)\\
%
&\in& (x,y) - \int_0^tKds +\Big(\int_0^tX_sr_sds, \int_0^tY_sb_sds\Big)+\Big(0, \int_0^tY_s \sigma_sdW_s\Big)\\
&=& (x,y)+ \int_0^t  [ (X_sr_s , Y_s b_s )-K  ]ds +\int_0^t (0, Y_s \sigma_s)dW_s.
\eeaa
Clearly, this is a set-valued SDI with coefficients:
\beaa
\left\{\ba{lll}
F(t, \o, \bx):= \mbox{diag}[r_t(\o), b_t(\o)]\bx -K;  \ms\\
g(t, \o, \bx):=\mbox{diag}[0, \si_t(\o)]\bx, 
\ea\right. \q (t, \o)\in[0,T], ~\bx=(x, y)\in\hR^2, 
\eeaa
where diag$[a, b]$ denotes the $2\times 2$ diagonal matrix with diagonal elements $a$ and $b$, and $C:=K$ is a given and fixed convex cone in $\hR^2$. Clearly, if the processes $r$, $b$ and $\si$ are $\hF$-adapted, then $F$ and $g$ satisfy Assumption \ref{assumpSDI}. A more general multi-dimnsional continuous-time super-hedging problem based on this model will be studied in our forthcoming work \cite{ACM}. 

\ms

\no {\bf Acknowledgement}: We would like to thank Prof. \c{C}a\u{g}{\i}n Ararat for many helpful discussions and suggestions on the subject during his visit to USC.

\end{document}